\documentclass[12pt]{amsart}
\usepackage{amssymb, amsmath, amscd,graphicx,subfigure}
\usepackage{epstopdf}

\newtheorem{thm}{Theorem}
\newtheorem{lem}[thm]{Lemma}
\newtheorem{cor}[thm]{Corollary}
\newtheorem{prop}[thm]{Proposition}

\theoremstyle{plain}

\theoremstyle{definition}

\newtheorem{defi}[thm]{Definition}
\newtheorem{construct}[thm]{Construction}

\theoremstyle{definition}
\newtheorem{rem}[thm]{Remark}

\renewcommand{\int}{\operatorname{int}}

\newcommand{\poff}{\operatorname{\mathbf {pull-off}}}
\newcommand{\pin}{\operatorname{\mathbf {push-in}}}

\newcommand{\lk}{\operatorname{lk}}

\newcommand{\Z}{\mathbb{Z}}
\newcommand{\G}{\mathbb{G}}
\newcommand{\W}{\mathbb{W}}

\newcommand{\R}{\mathbb{R}}

\newcommand{\Q}{\mathbb{Q}}

\newcommand{\cL}{\mathcal L}

\newcommand{\cW}{\mathcal{W}}
\newcommand{\cA}{\mathcal A}
\newcommand{\cB}{\mathcal B}

\newcommand{\e}{\epsilon}

\newcommand{\imra}{\looparrowright}

\newcommand{\ra}{\longrightarrow}

\begin{document}

\title[Jacobi Identities]{Jacobi Identities in Low-dimensional Topology}

%
%
\author{James Conant, Rob Schneiderman and Peter Teichner}
\address{Dept.~of Mathematics, University of Tennessee, Knoxville, TN 37996}
\email{jconant@math.utk.edu}
\address{Dept.~of Mathematics and Computer Science, Lehman College, City University of New York, Bronx, NY 10468 }
\email{robert.schneiderman@lehman.cuny.edu}
\address{Dept.~of Mathematics, University of California, Berkeley, CA 94720-3840}
\email{teichner@math.berkeley.edu}


\thanks{All authors are supported by the NSF.
This collaboration started during a joint visit at the Max-Planck
Institute in Bonn.}

\begin{abstract}
The Jacobi identity is the key relation in the definition of a Lie
algebra. In the last decade, it also appeared at the heart of the
theory of finite type invariants of knots, links and 3-manifolds
(and is there called the IHX relation). In addition, this relation
was recently found to arise naturally in a theory of embedding
obstructions for 2-spheres in $4$-manifolds in terms of Whitney
towers. This paper contains the first proof of the 4-dimensional
version of the Jacobi identity. We also expose the underlying
topological unity between the 3- and 4-dimensional IHX relations,
deriving from a beautiful picture of the Borromean rings embedded on
the boundary of an unknotted genus~3 handlebody in 3-space. This
picture is most naturally related to knot and 3-manifold invariants
via the theory of grope cobordisms.
\end{abstract}

\maketitle

\section{Introduction}

The only axiom in the definition of a Lie algebra, in addition to
the bilinearity and skew-symmetry of the Lie bracket, is the {\em
Jacobi identity}
\[
[[a,b],c]-[a,[b,c]]+[[c,a],b] =0.
\]
If the Lie algebra arises as the tangent space at the identity
element of a Lie group, the Jacobi identity follows from the
associativity of the group multiplication. Picturing the Lie
bracket as a rooted Y-tree with two inputs (the tips) and one
output (the root), the Jacobi identity can be encoded by the
following figure:

\begin{figure}[ht!]
         \centerline{\includegraphics[scale=.4]{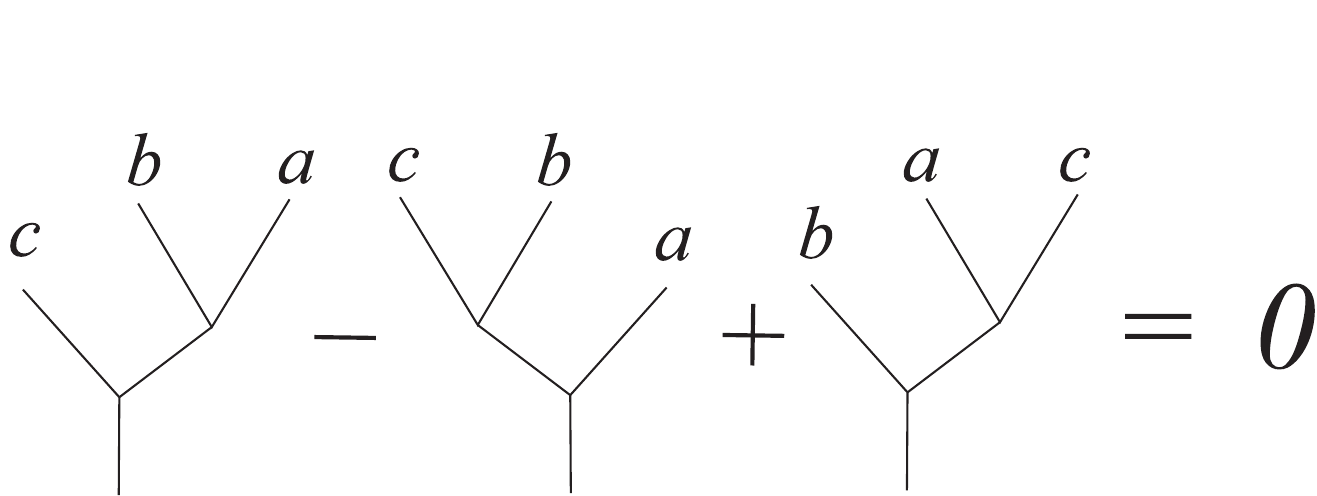}}
         \caption{The Jacobi identity}
         \label{Jacobi-fig}
\end{figure}

One should read this tree from top to bottom, and note that the
planarity of the tree (together with the counter-clockwise
orientation of the plane) induces an ordering of each trivalent
vertex which can thus be used as the Lie bracket. A change of this
ordering just introduces a sign due to the skew-symmetry of the
bracket. This will later correspond to the antisymmetry relation
for diagrams. Changing the input letters $a,b,c$ to $1,2,3$ and
labeling the root $4$, Figure~\ref{Jacobi-fig} may be redrawn with
the position of the labeled univalent vertices fixed as follows:

\begin{figure}[ht!]
         \centerline{\includegraphics[scale=.4]{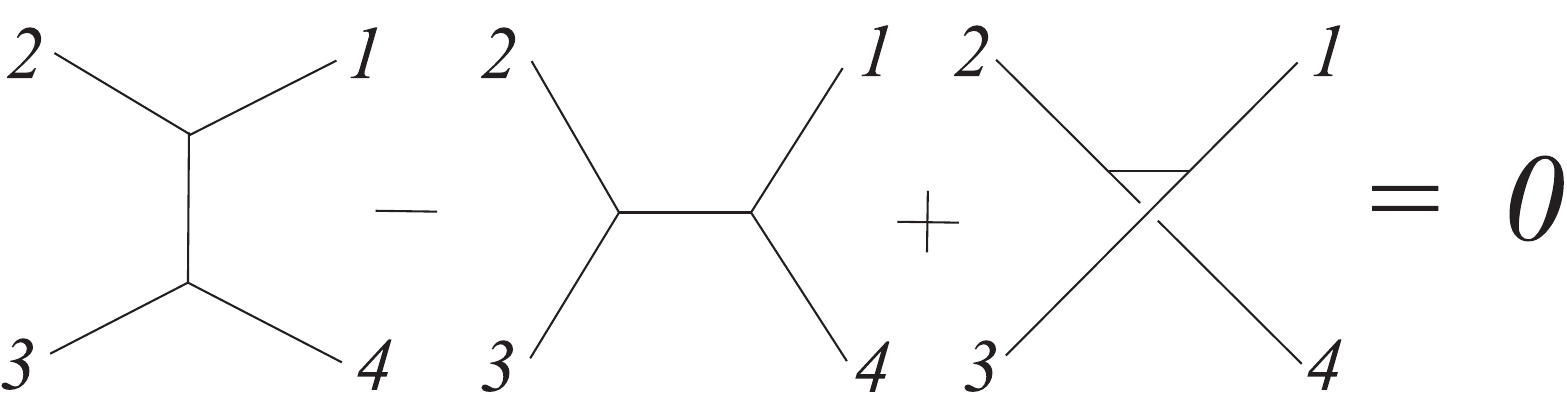}}
         \caption{The IHX-relation}
         \label{IHX-relation}
\end{figure}

This (local) relation is an unrooted version of the Jacobi
identity, and is well known in the theory of finite type (or
Vassiliev) invariants of knots, links and $3$-manifolds. Because
of its appearance it is called the IHX relation. The precise
connection between finite type invariants and Lie algebras is very
well explained in many references, see e.g. \cite{BN}.

Garoufalidis and Ohtsuki \cite{GO} were the first to prove a
version of a 3-dimensional IHX relation. It was needed to show
that
 a map from trivalent diagrams to homology 3-spheres was well-defined.
Habegger \cite{Ha} improved and conceptualized their construction.
Moving to the techniques of claspers (clovers), Garoufalidis,
Goussarov and Polyak \cite{GGP} sketch a proof of
Theorem~\ref{topihx} below, a theorem of which Habiro was also
aware. Our proof is completely new, and, we believe, more
conceptual. Moreover, it serves as a bridge between the 3- and
4-dimensional worlds.

\subsection{A Jacobi Identity in 4 Dimensions}
\begin{figure}[ht!]
         \centerline{\includegraphics[scale=.55]{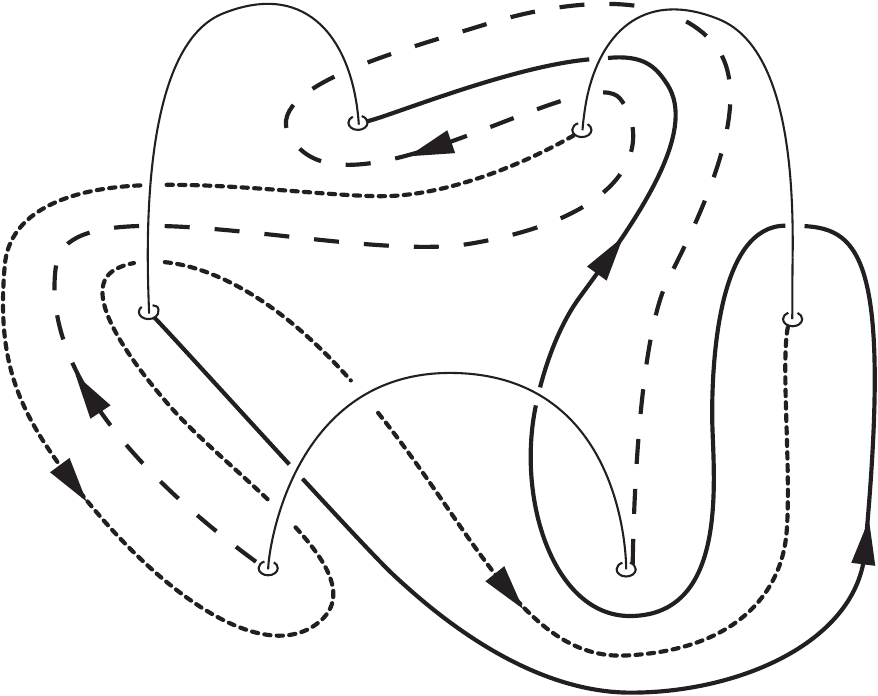}}
         \caption{The geometric origin of the Jacobi identity in Dimension 4.}
         \label{IHX-7B-fig}
\end{figure}

In Section~\ref{sec:4-dim} of this paper we will rediscover the
Jacobi relation in the context of intersection invariants for
Whitney towers in $4$-manifolds. It is actually a direct
consequence of a beautifully symmetric picture,
Figure~\ref{IHX-7B-fig}. The expert will see three standard
Whitney disks whose Whitney arcs are drawn in an unconventional
way (to be explained in Section \ref{sec:4d-IHX} below).
Ultimately, the freedom of choosing the Whitney arcs in this way
forces the Jacobi relation upon us.

The reader will recognize the 3-component link in the figure as
the Borromean rings. Each component consists of a semicircle and a
planar arc (solid, dashed, dotted respectively), exhibiting the
Borromean rings as embedded on the boundary of an unknotted
genus~3 handlebody in 3-space.

The Jacobi relation for Whitney towers plays a key role in the
obstruction theory for embedding $2$-spheres into $4$-manifolds
developed in \cite{ST1}. However, it was not proven in that
reference and the main purpose of this paper is to give a precise
formulation and proof of this Jacobi relation, see
Theorem~\ref{thm:4dihx} below. In sections \ref{sec:w-tower},
\ref{sec:int-trees} and \ref{sec:4d-IHX} of this paper, no
background is required of the reader beyond a willingness to try
to visualize surfaces in 4--space, and our elementary construction
can also serve as an introduction to Whitney towers.

Roughly speaking, a Whitney tower is a 2-complex in a
$4$-manifold, formed inductively by attaching layers of Whitney
disks to pairs of intersection points of previous surface stages,
see Section~\ref{sec:w-tower}. A Whitney tower has an {\em order}
which measures how many layers were used. Moreover, for any
unpaired intersection point $p$ in a Whitney tower $\cW$ of order
$n$, one can associate a tree $t(p)$ embedded in  $\cW$, see
Figure~\ref{IHX-W-tower-fig}. The tree $t(p)$ is a trivalent tree
with $n$ trivalent vertices, each representing a Whitney disk in
the tower. Each univalent vertex of $t(p)$ lies on a bottom stage
(immersed) sphere $A_i$ and is labeled by the index $i$.

Orientations of the surface stages in $\cW$ give
vertex-orientations of $t(p)$, i.e. cyclic orderings of the
trivalent vertices, and they also give a sign $\epsilon_p$.  We
define the {\em geometric intersection tree}
$\widetilde{\tau}_n(\cW)$ as the disjoint union of signed
vertex-oriented trees, one for each unpaired intersection point
$p$:
\[
\widetilde{\tau}_n(\cW):=\amalg_p \,\epsilon_p \cdot t(p).
\]
Properly interpreted, this union represents an obstruction to the
existence of an order $(n+1)$ Whitney tower extending $\cW$. (Note
that essentially the same geometric intersection tree is denoted
by $t_n(\cW)$ in \cite{ST1}.) The main result of this paper can
now be formulated as follows:

\begin{thm}(4-dimensional Jacobi relation) \label{thm:4dihx}
There exists an order 2 Whitney tower $\cW$ on four immersed
2-spheres in the 4--ball such that $\widetilde{\tau}_2(\cW)=(+I)
\amalg (-H) \amalg (+X)$ where $I$, $H$ and $X$ are the trees
shown in Figure~\ref{IHX-relation}.
\end{thm}

This result comes from the fact alluded to before, namely that
Whitney towers have the indeterminancy of choosing the Whitney
arcs! The local nature of Theorem~\ref{thm:4dihx} enables
geometric realizations of Jacobi relations via controlled
manipulations of Whitney towers, an essential step in the
obstruction theory described in \cite{ST1}. It should be mentioned
that there is also a 4-dimensional geometric Jacobi relation which
uses a Whitney move to locally replace an $I$-tree by an $H$-tree
and an $X$-tree (see \cite{S}).

In the easiest case $n=0$, a Whitney tower (of order~0) is just a
union of immersed $2$-spheres $A_1,\dots,A_\ell:S^2\imra M^4$, and
its geometric intersection tree $\widetilde{\tau}_0(\cup_i A_i)$
is a disjoint union of signed edges, one for each intersection
point among the $A_i$. The endpoints of the edges are labeled by
the $2$-spheres, or better by elements of the set
$\{1,\dots,\ell\}$, organizing the information as to which $A_i$
are involved in the intersection. Edges with index $i$ on both
ends correspond exactly to self-intersections of $A_i$.

In this case we know how to extract an invariant, namely by just
summing all the order 0 trees (= edges) in
$\widetilde{\tau}_0(\cup_i A_i)$, each signed edge of
$\widetilde{\tau}_0$ thought of as an integer $\pm 1$,  to get
exactly the {\em intersection numbers} among the $A_i$. Actually,
if $M$ is not simply connected, these ``numbers'' should be
evaluated in the group ring of $\pi_1M$, rather than in $\Z$,
leading to Wall's intersection invariants \cite{W}. This
corresponds to putting orientations and group elements on the
edges of each $t(p)$, and has been worked out in \cite{ST1} for
higher order Whitney towers. Note that for identical indices at
the ends of an edge, the two possible orientations on the edge
give isomorphic pictures, leading to the usual relations in the
group ring when measuring self-intersections:
\[
w_1(g)\cdot g = g^{-1}, \quad \forall g\in \pi_1M.
\]

In the present paper our constructions are local so that we may
safely ignore these group elements.

If $\widetilde{\tau}_0(\cup_i A_i)$ sums to zero then all the
intersections can be paired up by Whitney disks, i.e. there is a
Whitney tower $\cW$ of order~1 with the $A_i$ as bottom stages.
Then $\widetilde{\tau}_1(\cW)$ is a disjoint union of signed
(vertex-oriented) Y-trees, and again the univalent vertices have
labels from $\{1,\dots,\ell\}$. It was shown in \cite{ST} (and in
\cite{Ma}, \cite{Y} for simply connected $4$-manifolds) that a
summation as above leads to an invariant $\tau_1(\cW)$ which
vanishes if and only if there is a Whitney tower $\cW$ of order~2
with the $A_i$ as bottom stages. In fact, if defined in the
correct target group, $\tau_1(\cW)$ only depends on the regular
homotopy classes of the $A_i$ and hence is a well defined higher
obstruction for representing these classes by disjoint embeddings.

Theorem~\ref{thm:4dihx} only becomes relevant for
$\widetilde{\tau}_2$ and higher, and we next give a proper
formulation of some necessary notation and terminology for
intersection trees.

\begin{defi} \label{def:degree}
We define the {\em order} of a trivalent tree to be the number of
trivalent vertices and the {\em degree} to be one \emph{more} than
that number. The degree is also one half of the total number of
vertices, or one less than the number of univalent vertices. This
definition is consistent with the theory of finite type
invariants, where the degree goes back to Vassiliev.
\end{defi}

\begin{defi} \label{def:B}
Consider pairs $(\epsilon,t)$ where $\epsilon\in\{\pm\}$ and $t$
is a vertex-oriented trivalent tree of degree $n$, with univalent
vertices labeled from $\{1,\dots,\ell\}$.
\begin{enumerate}
\item An AS (antisymmetry)
relation is of the form
\[
(\epsilon,t)=(-\epsilon,t'),
\]
where $t'$ is isomorphic to $t$ and its orientation differs from
that of $t$ by changing the cyclic orientation at a single vertex.
All AS relations generate an equivalence relation, and we let
$\widetilde{\cB}^t_{n}(\ell)$ be the commutative monoid with unit
generated by the set of equivalence classes of such pairs
$(\epsilon,t)$. We think of the monoid operation as disjoint
union, $\amalg$, and we write $\epsilon\cdot t$ for the
equivalence class of $(\epsilon,t)$.

\item The abelian group $\widehat{\cB}^t_n(\ell)$  is obtained by dividing
the monoid $\widetilde{\cB}^t_n(\ell)$ by all relations of the
form
\[
(\epsilon\cdot t) \amalg (-\epsilon \cdot t)=0,
\]
where $0$ is the unit of the monoid $\widetilde{\cB}^t_n(\ell)$.
This clearly introduces inverses and the monoid operation $\amalg$
becomes a group addition which we write as ``+''.
\item The abelian group $\cB^t_n(\ell)$ is obtained
from $\widehat{\cB}^t_n(\ell)$ by dividing out all Jacobi (IHX)
relations.
\end{enumerate}
\end{defi}

\begin{rem} \label{rem:concrete}
Definition (i) can be spelled out more concretely at two points:
The equivalence relation {\em generated by} AS relations as above
is just given by relations of the form
\[
(\epsilon,t)=((-1)^k\epsilon,t')
\]
where $t'$ is isomorphic to $t$ and its orientation differs from
that of $t$ by changing the cyclic orientation at exactly k
vertices. Moreover, the commutative monoid
$\widetilde{\cB}^t_{n}(\ell)$ {\em generated by} such equivalence
classes can be described by working with (equivalence classes of)
{\em finite unions} of vertex-oriented trees, with each connected
component labeled by a sign $\epsilon$. Then the disjoint union
really gives a monoid structure on this set which is clearly
commutative and generated by trees. Its unit is given by the empty
graph.
\end{rem}

Let $\W_{(n-1)}(\ell)$ denote the set of Whitney towers of order
$(n-1)$ on bottom stages $A_1,\dots,A_\ell$. We have been
discussing a map $\widetilde{\tau}_{n-1}$ which we can now write
as
$$\widetilde{\tau}_{n-1}\colon \W_{(n-1)}(\ell)\to \widetilde{\cB}^t_n(\ell).
$$
Working modulo the relations in definition (ii) above, we get a
summation map
\[
\widehat\tau_{(n-1)}\colon \W_{(n-1)}(\ell)\ra
\widehat{\cB}_n^t(\ell).
\]
More explicitly, if $\widetilde{\tau}_{n-1}(\cW)=\amalg_p \,\e_p
\cdot t(p)$ is the geometric intersection tree of an order $(n-1)$
Whitney tower $\cW$ we set
$$
\widehat{\tau}_{(n-1)}(\cW):=\sum_p \e _p \cdot
t(p)\in\widehat{\cB}_n^t(\ell)
$$

It is a consequence of the AS relation that only orientations of
the {\em bottom stages} $A_i$ are relevant for the definition of
$\widetilde\tau_{(n-1)}(\cW)$, see Lemma~\ref{lem:W-tower-AS}.
From our geometric point of view, this is the main reason for
introducing AS relations.

 The question arises as to
whether $\widehat\tau_{(n-1)}(\cW)$ can be made into an
obstruction for representing the bottom stages $A_i$, up to
homotopy, by disjoint embeddings. The punch-line of the first part
of this paper is that this can only be possible if we quotient the
groups  $\widehat{\cB}^t_n(\ell)$ by all Jacobi relations. This
gives the above-mentioned groups $\cB^t_n(\ell)$, containing
elements $\tau_{(n-1)}(\cW)$, which are more customary in the
theory of finite type invariants. In fact, in the general finite
type theory the superscript `$t$', for tree, does not appear
because one uses all trivalent graphs instead of just unions of
trees.

Theorem~\ref{thm:4dihx} implies that one needs to study these
quotients if one wants to obtain invariants of the bottom spheres
$A_i$ from the intersection trees associated to Whitney towers. As
shown in \cite{ST1}, the vanishing of $\tau_{(n-1)}(\cW)$ in
$\cB^t_n(\ell)$ is sufficient for finding a next order $n$ Whitney
tower on the $A_i$ (up to homotopy). However, it is an open
problem what precise further quotient of $\cB^t_n(\ell)$ is
necessary to get a well-defined invariant which only depends on
the homotopy classes of the $A_i$ (and {\em not} on the
order~$(n-1)$ Whitney tower) and gives the complete obstruction to
the existence of an order~$n$ Whitney tower.

 \subsection{From 4- to 3-dimensional Jacobi relations}\label{sec:down-to-3}
In Section~\ref{sec:3and4} we connect the geometric Jacobi
relation explained above to a 3-dimensional setting via a
correspondence between capped grope concordances and Whitney
towers. This translation becomes important because, up to date,
there is no useful definition of a Whitney tower in 3 dimensions.
On the other hand, two of us have introduced in \cite{CT1} a
theory of (capped) grope cobordisms between knots in 3-space, and
the third member in our group  \cite{S} has worked out a
4-dimensional correspondence between capped grope concordances and
Whitney towers.

A grope is a certain 2-complex, built out of layers of surfaces.
The number of these layers is measured by the {\em class} of the
grope, later corresponding to the degree of a tree. A grope
contains a specified {\em bottom stage surface}, usually with one
or two boundary circles, depending on whether it is used to relate
string links or links. This is explained in detail in
Section~\ref{sec:gropes} and we shall introduce the notation that
\begin{itemize}
\item a grope {\em cobordism} is an embedding of a grope into 3-space, see Section~\ref{sec:grope cobordism}.
\item a grope {\em concordance} is an embedding of a grope  into 4-space.
More precisely, the embedding is into $B^3 \times [0,1]$, see
Section~\ref{sec:3to4}.
\end{itemize}
We shall also explain the notions of {\em capped} grope cobordism
and concordance. The {\em caps} are embedded disks whose
boundaries lie on the top stages of the grope. The (interiors of
the) caps are allowed to intersect the grope only in the bottom
stage surface. The punch-line is that these intersections are
going to be
\begin{itemize}
\item arcs, from one part of the boundary to another, in the bottom stage of a grope cobordism,
\item points in the bottom stage of a grope concordance.
\end{itemize}
These statements are the generic case in dimension~4 and need
certain cleaning up operations in dimension~3. In any case, when
pushing a grope cobordism into 4-space, the arcs become points,
and one loses the information of the order in which the arcs hit
the boundary. More precisely, in Section~\ref{sec:3to4} we shall
explain in full detail the following commutative diagram:

\begin{equation}\label{diagram}
\begin{CD}
 \G^c_n(\ell)  @>{\pin}>>  \overline{\W}_{(n-1)}(\ell)  \\
  @V{\widetilde\tau^c_n}VV @V{\widetilde{\tau}_{(n-1)}}VV \\
  \widetilde{\cA}^t_n(\ell)  @>{\poff}>> \widetilde{\cB}^t_n(\ell)
\end{CD}
\end{equation}

Here $\G^c_n(\ell)$ is the set of class $n$ capped grope
cobordisms of $\ell$-string links. The set $
\overline{\W}_{(n-1)}(\ell)$ is a quotient of $\W_{(n-1)}(\ell)$
by the relation equating Whitney towers which are assigned the
same element by $\widetilde{\tau}_{(n-1)}$. The map $\pin$ takes a
capped grope, pushes it slightly into the 4-ball, and then surgers
the resulting grope concordance into a Whitney tower (of order
$(n-1)$). This procedure has some non-uniqueness which is why we
need the space  $ \overline{\W}_{(n-1)}(\ell)$ as opposed to
$\W_{(n-1)}(\ell)$. The monoid $\widetilde{\cA}^t_n(\ell)$ is just
like its $\cB$-analogue, except that the univalent vertices of the
trees are attached to $\ell$ numbered strands (which form a
trivial string link). More precisely, we have

\begin{defi} \label{def:A}
Consider pairs $(\epsilon,t)$ where $\epsilon\in\{\pm\}$, and $t$
is a vertex-oriented trivalent tree of degree $n$ whose tips are
attached to the trivial $\ell$-string link.
\begin{enumerate}
\item  As in Definition~\ref{def:B}, AS (antisymmetry)
relations of these pairs generate an equivalence relation, and we
let $\widetilde{\cA}^t_n(\ell)$ be the abelian monoid generated by
the equivalence classes. As before, the monoid operation is given
by disjoint union and we write $\epsilon\cdot t$ for
$(\epsilon,t)$.
\item The abelian group $\widehat{\cA}^t_n(\ell)$ is obtained from
$\widetilde{\cA}^t_n(\ell)$ by dividing by all relations of the
form
\[
(\epsilon\cdot t)\amalg (-\epsilon\cdot t)=0,
\]
where $0$ is the trivial $\ell$-string link with the `empty graph'
attached. The monoid operation $\amalg$ becomes the group addition
``$+$''.
\item The abelian group $\cA^t_n(\ell)$ is obtained from
$\widehat{\cA}^t_n(\ell)$ by dividing out all Jacobi (IHX)
relations.
\end{enumerate}
\end{defi}

The homomorphism $\poff$ in diagram~(\ref{diagram}) pulls each
tree off of the strands of the trivial $\ell$-string link, just
remembering their indices in $\{1,\dots,\ell\}$. Thus the diagram
above says exactly what information is lost when one moves from 3
to 4 dimensions, namely the orders in which the  caps hit the
bottom stages.

The map $\widetilde\tau^c_n$ is defined precisely in
Definition~\ref{def:cappedtau} using a notion of geometric
intersection trees for gropes (Definition~\ref{grope tree}), and
just as in the case of Whitney towers, leads to maps
$\widehat\tau^c_n$ and $\tau^c_n$.

By re-interpreting our central picture, Figure~\ref{IHX-7B-fig},
in terms of capped gropes in 3-space, $\widetilde\tau^c_n$ will be
used to show that the 4-dimensional Jacobi relation from
Theorem~\ref{thm:4dihx} can be lifted to a 3-dimensional version:

\begin{thm} (3-dimensional Jacobi relation)\label{thm:3dihx}
Suppose $t_I,t_H,t_X$ are the three terms in any IHX relation
 in $\widetilde{\cA}^t_3(\ell)$. Then there is a
class~3 simple grope cobordism $G^c$, which takes the
$\ell$-component trivial string link to itself, such that
$\widetilde\tau^c_3(G^c)=(+t_I)\amalg(-t_H)\amalg(+t_X)$.
\end{thm}

We should remark that by a main theorem of \cite{CT1}, we can
think of $\G^c_n(\ell)$ as being the set of degree $n$ capped (or
simple) claspers in the complement of some $\ell$ component
(string) link. The map $\widetilde\tau^c_n$ is then the obvious
map which sends a clasper to its tree type, with univalent
vertices attached to the link components which link with the
corresponding tips. However, $\widetilde\tau^c_n$ can also be
directly defined for capped gropes, as we explain in
Definition~\ref{def:cappedtau}. One consequence of our work is the
following theorem:

\begin{thm}(3-dimensional Jacobi relation for claspers)\label{topihx}
Suppose three tree claspers $C_i$ differ locally by the three
terms in the Jacobi relation. Given an embedding of $C_1$ into a
$3$-manifold, there are embeddings of $C_2$ and $C_3$ inside a
regular neighborhood of $C_1$, such that the leaves are parallel
copies of the leaves of $C_1$, and the edges avoid any caps that
$C_1$ may have. Moreover, surgery on $C_1\cup C_2\cup C_3$ is
diffeomorphic (rel boundary of the handlebody neighborhood) to
doing no surgery at all.
\end{thm}

This theorem was stated and utilized in \cite{CT2}, although the
fact that the claspers must be tree claspers was accidentally
omitted. The theorem was needed  to prove Lemma 3.11(a) in
\cite{CT2}. We reprove this lemma as Lemma~\ref{ct2lem} of the
current paper, as an easy consequence of our general machinery.
The map $\tau_n^c$ is relevant to the theory of Vassiliev
invariants. Given a simple grope cobordism between two links, it
records the difference in the Vassiliev invariants of degree $n$
between the two links. Thus, similar in spirit to $\tau_n$ for
Whitney towers, $\tau_n^c$ could represent an obstruction to two
links being isotopic. However, again the question remains on how
it depends on the particular choice of the given grope cobordism.

\subsection{Grope cobordism of string links}
In the last Section~\ref{sec:string}, we shall use the techniques
developed in this paper to obtain new information about string
links. Let $\cL(\ell)$ be the set of isotopy classes of string
links in $D^3$ with $\ell$ components (which is a monoid with
respect to the usual ``stacking'' operation). Its quotient by the
relation of grope cobordism of class $n$ is denoted
$\cL(\ell)/G_n$, compare Definition~\ref{grope cobordism}. The
quotient by the relation of \emph{capped} grope cobordism of class
$n$ is denoted by $\cL(\ell)/G^c_n$. The submonoid of $\cL(\ell)$,
consisting of those string links which cobound a class $n$  grope
with the trivial string link, is denoted by $G_n(\ell)$, and
similarly the submonoid consisting of those string links which
cobound a class $n$  capped grope with the trivial string link, is
denoted by $G^c_n(\ell)$. The relation of capped (respectively not
capped) grope cobordism of class $n$ coincides with the relation
that two string links differ by a sequence of simple (respectively
rooted) clasper surgeries of degree $n$. Using this connection and
results of Habiro \cite{H} we show

\begin{thm}\label{thm:nil}
$\cL(\ell)/G_{n+1}$ and $\cL(\ell)/G^{c}_{n+1}$  are finitely
generated groups and the iterated quotients
\[
G_n(\ell)/G_{n+1} \quad \text{ respectively } \quad
G^{c}_n(\ell)/G^{c}_{n+1}
\]
are central subgroups. As a consequence,  $\cL(\ell)/G_{n+1}$ and
$\cL(\ell)/G^{c}_{n+1}$ are nilpotent.
\end{thm}

In the case of knots, $\ell=1$, results of Habiro and also
\cite{CT2} imply that $G^c_n(1)/G^c_{n+1}$ is rationally
isomorphic to the space $\mathcal B_n\otimes \Q$ appearing in the
theory of Vassiliev invariants (Indeed, we alluded to $\mathcal
B_n = \mathcal B_n(1)$
 a few paragraphs after Definition~\ref{def:B}).

For the case of $\ell\geq 2$ no such theorem is known, but we show
that if one relaxes the requirement that the gropes be capped
(which is the same as relaxing the requirement that all leaves of
the clasper bound disjoint disks to the requirement that only one
leaf does) then one does get such a statement. Using our geometric
IHX relations, we will construct a surjective homomorphism from
diagrams to string links modulo grope cobordism:
$$
\Phi_n\colon \cB^g_n(\ell)\twoheadrightarrow G_n(\ell)/G_{n+1}
$$
where $\cB^g_n$ denotes the usual abelian group of trivalent
graphs, modulo IHX- and AS-relations, but graded by the {\em grope
degree} (which is the Vassiliev degree plus the first Betti number
of the graph), compare Section~\ref{sec:genihx}.

\begin{thm}\label{kontsevich-grope}
The map
$\Phi_n\otimes\Q\colon\cB^g_n(\ell)\otimes\Q\overset{\cong}{\ra}
G_n(\ell)/G_{n+1}\otimes\Q$ is an isomorphism.
\end{thm}

 This extends a result in
\cite{CT2} from knots to string links and it relies on the
existence of the Kontsevich integral for string links, which
serves as an inverse to the above map. Although Theorem
\ref{kontsevich-grope} is an elementary modification of the
argument in \cite{CT2}, we found it to be quite surprising in
light of the fact that the corresponding statement for capped
gropes and simple claspers is unknown.

The map $\Phi_n$ comes from a map $\widehat{\tau}^g_n$ defined in
Section~\ref{sec:genihx}, which assigns a linear combination of
vertex-oriented unitrivalent graphs of grope degree $n$ to any
grope cobordism of class $n$. This map is a technical improvement
of our methods in \cite{CT2}, and is necessary for us to realize
the IHX relation in the uncapped case. To define the map in that
paper, we first turned a grope cobordism into a sequence of simple
clasper surgeries, and then read off the unitrivalent graphs from
the graph types of the claspers. In this paper, we read off the
graphs directly from the (genus one) grope itself.
 The proof that this map induces an isomorphism still
requires the techniques of \cite{CT2}, and in particular, still
requires the passage to claspers, since the Kontsevich integral's
behavior with respect to claspers is well understood.

In an appendix we define the map $\widehat{\tau}^g_n$ for
arbitrary grope cobordisms, which is more general than the genus
one gropes used in the body of the paper. This is logically not
necessary but included for completeness and possibly for future
use.

\vspace{1em} \noindent {\bf Acknowledgment}: It is a pleasure to
thank Tara Brendle, Stavros Garoufalidis and the referee for
helpful discussions.

\section{A Jacobi Identity in Dimension 4}\label{sec:4-dim}
In this section we prove Theorem~\ref{thm:4dihx}, but we first
explain some background material and state an important Corollary
which is used in \cite{ST1}. For more details on immersed surfaces
in 4--manifolds we refer to \cite{FQ}, for more details on Whitney
towers compare \cite{S}, \cite{S2}, \cite{ST1}.

\subsection{Whitney towers}\label{sec:w-tower}
\begin{figure}[ht!]
         \centerline{\includegraphics[scale=.65]{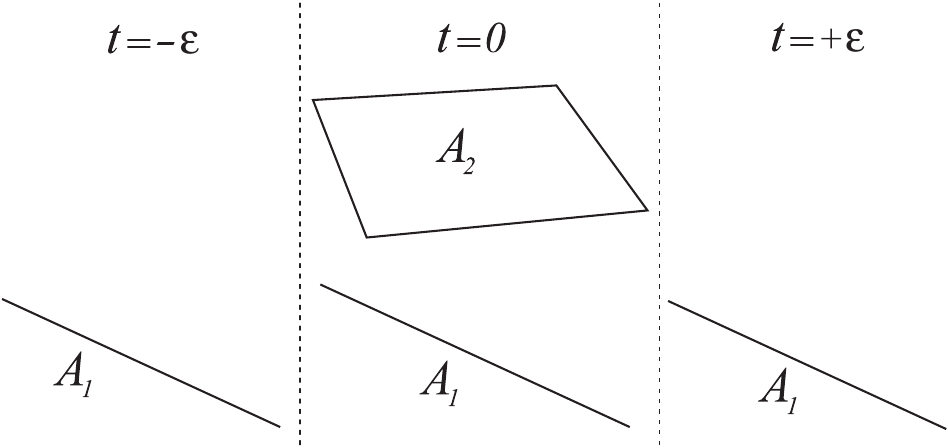}}
         \caption{}
         \label{IHX-Fingermove1-fig}
\end{figure}

Using local coordinates $\mathbb{R}^3\times(-\epsilon,+\epsilon)$,
Figure~\ref{IHX-Fingermove1-fig} shows a pair of disjoint local
sheets of oriented surfaces $A_1$ and $A_2$ in 4--space. We think
of the fourth coordinate as ``time'', so the sheet $A_2$ lies
completely in the \emph{present} $t=0$, whereas $A_1$ moves
through time and thus also forms a 2-dimensional sheet represented
by an arc which extends from \emph{past} into \emph{future}.
\begin{figure}[ht!]
         \centerline{\includegraphics[scale=.65]{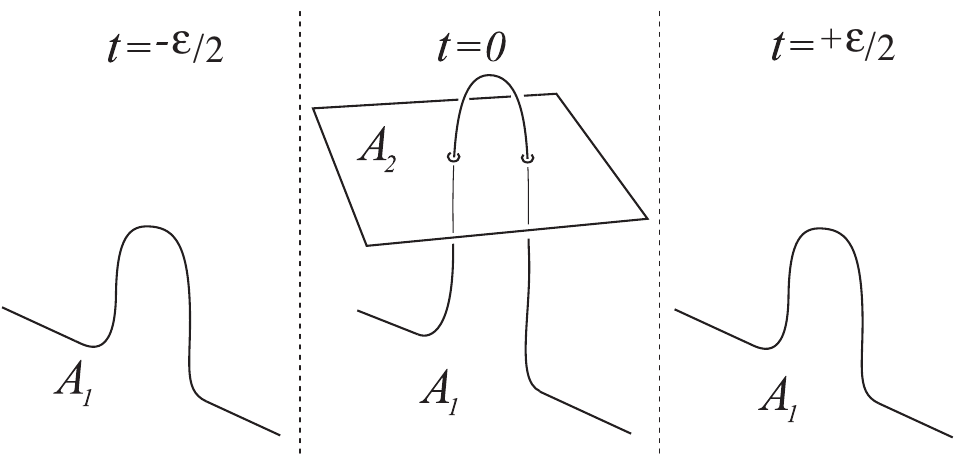}}
         \caption{}
         \label{IHX-Fingermove2-fig}
\end{figure}
\begin{figure}[ht!]
         \centerline{\includegraphics[scale=.6]{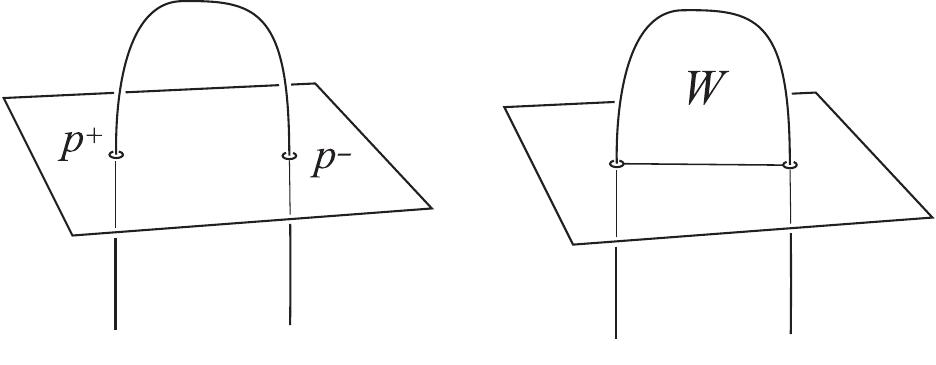}}
         \caption{Left: A cancelling pair of intersections $p^{\pm}$. Right: A
         Whitney disk pairing $p^{\pm}$.}
         \label{IHX-int-pair-W-disk-fig}
\end{figure}
Figure~\ref{IHX-Fingermove2-fig} shows the result of applying a
(Casson) {\em finger move} to the sheets of
Figure~\ref{IHX-Fingermove1-fig}, with $A_1$ having been changed
by an isotopy supported near an arc from $A_1$ to $A_2$, creating
a pair of transverse intersection points in $A_1\cap A_2\subset
\mathbb{R}^3\times\{0\}$. Such a pair of intersection points is
called a {\em cancelling pair} since their signs differ and they
can be paired by a {\em Whitney disk} as illustrated in
Figure~\ref{IHX-int-pair-W-disk-fig}. Note that the boundary of
the Whitney disk is a pair of arcs, one in each sheet, connecting
the cancelling pair of intersections. A Whitney disk guides a
motion (of either sheet) called a {\em Whitney move} that
eliminates the pair of intersection points \cite{FQ}. A Whitney
move guided by a Whitney disk whose interior is free of
singularities can be thought of as an ``inverse'' to the finger
move since it eliminates a cancelling pair without creating any
new intersections. In general, Whitney disks may have interior
self-intersections and intersections with other surfaces so that
eliminating a cancelling pair via a Whitney move may also create
new singularities. Pairing up ``higher order'' interior
intersections in a Whitney disk by ``higher order'' Whitney disks
leads to the notion of a Whitney tower:

\begin{defi}[(compare \cite{S},\cite{S2},\cite{ST1})]\label{defi:w-tower}\mbox{}

\begin{itemize}
\item A {\em surface of order 0} in a 4--manifold $M$
is an oriented surface in $M$ with boundary embedded in the
boundary and interior {\em immersed} in the interior of $M$. A
{\em Whitney tower of order 0}  is a collection of order 0
surfaces. These are usually referred to as the {\em bottom stage
surfaces} or \emph{underlying surfaces}, and a (higher order)
Whitney tower is built {\em on} these surfaces.
\item The {\em order of a (transverse) intersection point} between a surface of order $n$ and a
surface of order $m$ is $n+m$.
\item The {\em order of a Whitney disk} is $(n+1)$ if it pairs intersection points of order $n$.
\item For $n\geq 1$, a {\em Whitney tower of order $n$}  is a Whitney tower $\cW$ of
order~$(n-1)$ together with order~$n$ Whitney disks pairing all
order~$(n-1)$ intersection points of $\cW$, see
Figure~\ref{IHX-W-tower-fig}. These order~$n$ Whitney disks are
allowed to self-intersect, and/or intersect each other, as well as
lower order surfaces.
\end{itemize}
The boundaries of the Whitney disks in a Whitney tower are
required to be disjointly embedded and the Whitney disks
themselves are required to be {\em framed}. \hfill$\Box$
\end{defi}

Framings of Whitney disks will not be discussed here, see e.g.
\cite{FQ}. In the construction of an order 2 Whitney tower (proof
of Theorem \ref{thm:4dihx}) the reader familiar with framings can
check that the Whitney disks are correctly framed.

\begin{figure}[ht!]
         \centerline{\includegraphics[scale=.5]{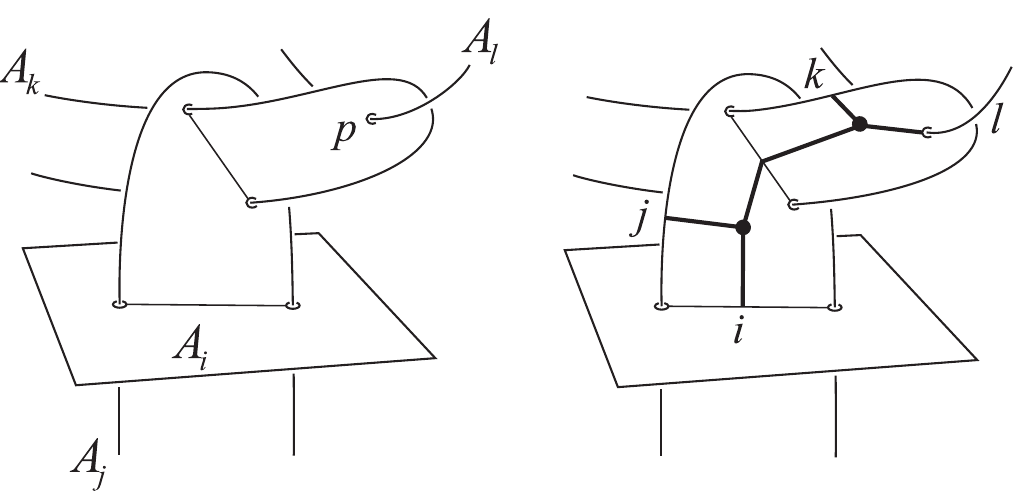}}
         \caption{Part of an order~2 Whitney tower on order~0 surfaces $A_i$, $A_j$, $A_k$, and
         $A_l$,
         and the labeled  tree $t(p)$ of order 2 = Vassiliev degree 3, associated to the
         order~2
         intersection point $p$.}
         \label{IHX-W-tower-fig}
\end{figure}

\subsection{Intersection trees for Whitney towers}\label{sec:int-trees}
For each order $n$ intersection point $p$ in a Whitney tower $\cW$
there is an associated labeled trivalent tree $t(p)$ of order $n$
(Figure~\ref{IHX-W-tower-fig}). The order of a tree is the number
of trivalent vertices (which is one less than the Vassiliev
degree).
 This tree $t(p)$ is most easily
described as a subset of $\cW$ which ``branches down'' from $p$ to
the order 0 surfaces, bifurcating in each Whitney disk: The
trivalent vertices of $t(p)$ correspond to Whitney disks in $\cW$,
the labeled  univalent vertices of $t(p)$ correspond to the
labeled  order 0 surfaces of $\cW$ and the edges of $t(p)$
correspond to sheet-changing paths between adjacent surfaces in
$\cW$.

Fixing orientations on the surfaces in $\cW$ (including Whitney
disks) endows each intersection point $p$ with a sign
$\e_p\in\{\pm\}$, determined as usual by comparing the
orientations of the intersecting sheets at $p$ with that of the
ambient manifold. These orientations also determine a cyclic
orientation for each of the trivalent vertices of $t(p)$ via a
bracketing convention which will be illustrated explicitly during
the proof of Theorem~\ref{thm:4dihx} below. We shall henceforth
assume that our Whitney towers come equipped with such
orientations.

The order~$n$ intersection points are the ``interesting''
intersection points in an order~$n$ Whitney tower $\cW$, since
these points may represent an obstruction to the existence of an
order~$(n+1)$ Whitney tower on the order~0 surfaces of $\cW$. (In
fact, all intersections of order greater than $n$ can be
eliminated by finger moves on the Whitney disks.)

Recall $\widetilde{\cB}^t_{n+1}$ from Definition~\ref{def:B}.
\begin{defi}\label{defi:t(W)}
 For an oriented order $n$ Whitney tower $\cW$, define $\widetilde{\tau}_n(\cW)\in\widetilde{\cB}^t_{n+1}(\ell)$,
the order $n$ {\em geometric intersection tree} of $\cW$, to be
the disjoint union of signed labeled  vertex-oriented trivalent
trees
$$
\widetilde{\tau}_n(\cW):=\amalg_p \,\e _p \cdot t(p)
$$
over all order~$n$ intersection points $p\in \cW$.\hfill$\Box$
\end{defi}
We emphasize that $\widetilde{\tau}_n(\cW)$ is a collection of
signed trees of order~$n$, possibly with repetitions, {\em
without} cancellation of terms. (The geometric intersection tree
is denoted by $t_n(\cW)$ in \cite{ST1}, as is an un-oriented
version in \cite{S}.)

Note that there is a natural map
$\pi\colon\widetilde{\cB}^t_{n+1}(\ell)\to
\widehat{\cB}^t_{n+1}(\ell)$ given by sending the monoid operation
to the group addition.

\begin{defi}
Given an oriented order $n$ Whitney tower $\cW$, define
$\widehat{\tau}_n(\cW)=\pi(\widetilde{\tau}_n(\cW))$.
\end{defi}

It turns out (see Lemma~\ref{lem:W-tower-AS} below) that for any
fixed Whitney tower $\cW$, the AS antisymmetry relations
correspond {\em exactly} to the indeterminacies coming from
orientation choices on the Whitney disks in $\cW$, so that the
element $\widetilde{\tau}_n(\cW)\in\widetilde{\cB}^t_{n+1}$  only
depends on the orientations of the bottom stage surfaces. On the
other hand, by fixing the bottom stage surfaces and varying the
choices of Whitney disks we are led to the IHX relations, as we
describe in the next section.

Since the ultimate goal of studying Whitney towers is to extract
homotopy invariants $\tau_n$ of the underlying order zero surfaces
from the geometric intersection tree, such an element should
vanish for any Whitney tower $\cW$ on immersed 2--spheres into
4--space since all such spheres are null-homotopic.
Theorem~\ref{thm:4dihx} from the introduction (proven below), and
its corollary (Corollary~\ref{4dihx-cor}) illustrate the necessity
of the IHX relation in the target of $\tau_n$. Since
Theorem~\ref{thm:4dihx} is a local statement (taking place in a
4--ball) it can be used to ``geometrically realize'' all higher
degree IHX relations for Whitney towers in arbitrary 4--manifolds,
a key part of the obstruction theory described in \cite{ST1}. The
following corollary of Theorem~\ref{thm:4dihx} is proved in
\cite{ST1}.

\begin{cor}\label{4dihx-cor}
Let $\cW$ be an order~$n$ Whitney tower on surfaces $A_i$. Then,
given any order~$n$ trivalent trees $t_I$, $t_H$ and $t_X$
differing only by a local IHX relation, there exists an order~$n$
Whitney tower $\cW'$ on $A'_i$ homotopic (rel boundary) to the
$A_i$ such that
$$
\widetilde{\tau}_n(\cW')=\widetilde{\tau}_n(\cW) \amalg (+t_I)
\amalg (-t_H) \amalg (+t_X).
$$
$\hfill\square$
\end{cor}

The idea of the proof of Corollary~\ref{4dihx-cor} is that by
applying finger moves to surfaces in a Whitney tower one can
create clean Whitney disks which are then tubed into the spheres
in Theorem~\ref{thm:4dihx}. This construction can be done without
creating extra intersections since finger moves are supported near
arcs and the construction of Theorem~\ref{thm:4dihx} is contained
in a 4--ball.

\subsection{Proof of the Main Theorem~\ref{thm:4dihx}}\label{sec:4d-IHX}

\begin{figure}[ht!]
         \centerline{\includegraphics[scale=.55]{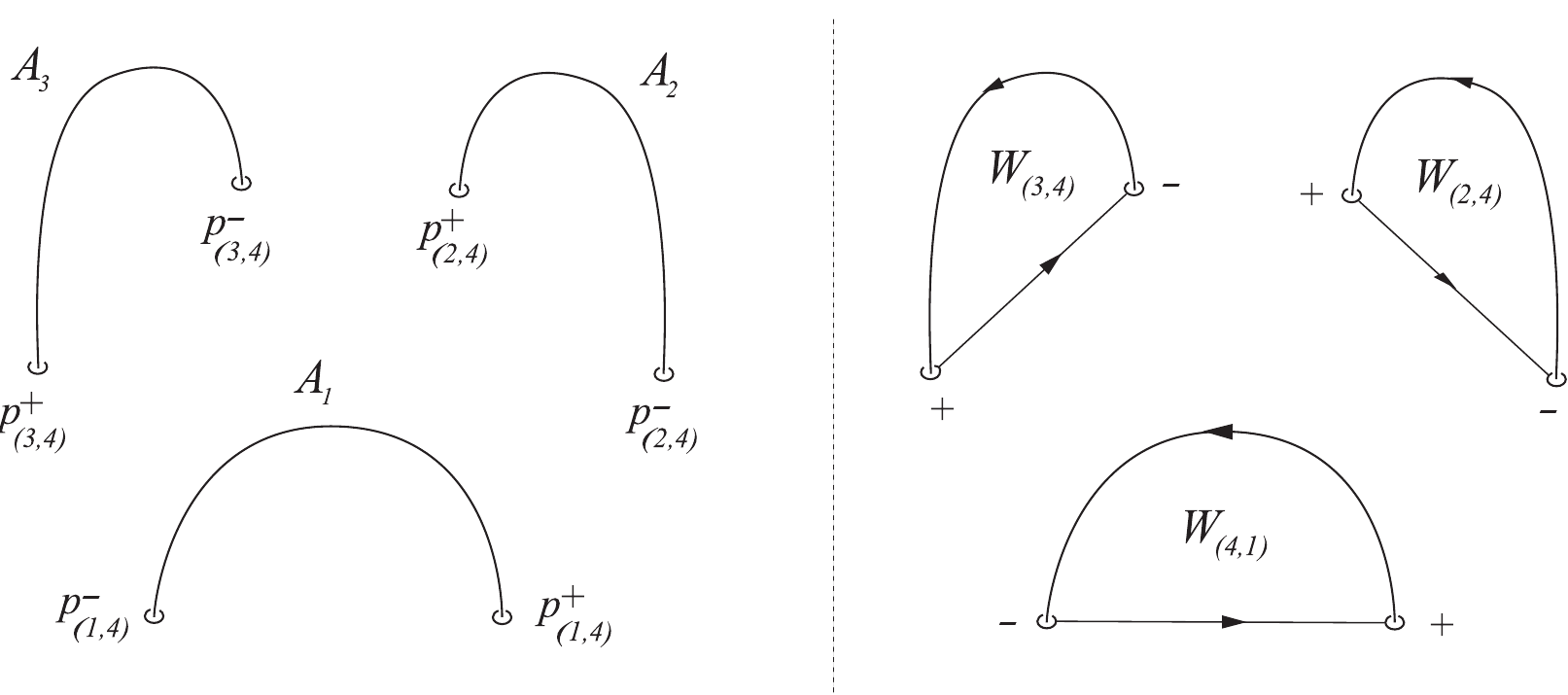}}
         \caption{The clean order~1 Whitney tower $\cW^0$ is shown on the right.}
         \label{IHX-1+3-fig}
\end{figure}

The 4--dimensional IHX construction starts with any four
disjointly embedded oriented 2-spheres $A_1,A_2,A_3,A_4$ in
4--space. Perform finger moves on each $A_i$, for $i=1,2,3$, to
create a cancelling pair of order zero intersection points
$p^{\pm}_{(i,4)}$ between each of the first three 2-spheres (still
denoted $A_i$) and $A_4$ as pictured in the left-hand side of
Figure~\ref{IHX-1+3-fig} where $A_4$ appears as the ``plane of the
paper'' with the standard counter-clockwise orientation, sitting
in the ``present'' slice $\R^3\times\{0\}$ of local coordinates
$\R^3\times (-\epsilon,+\epsilon)$ in 4-space. Choose disjointly
embedded oriented order~1 Whitney disks $W_{(3,4)}$, $W_{(2,4)}$
and $W_{(4,1)}$ for the cancelling pairs $p^{\pm}_{(i,4)}$ as in
the right-hand side of Figure~\ref{IHX-1+3-fig}. Here the bracket
sub-script notation corresponds to the following {\em orientation
convention}: The bracket subscript $(i,j)$ on a Whitney disk
indicates that the boundary $\partial W_{(i,j)}$ of the Whitney
disk is oriented from the negative intersection point to the
positive intersection point along $A_i$ and from the positive to
the negative intersection point along $A_j$. This orientation of
$\partial W_{(i,j)}$ together with a second ``inward pointing''
tangent vector induces the orientation of $W_{(i,j)}$.
\begin{figure}[ht!]
         \centerline{\includegraphics[scale=.55]{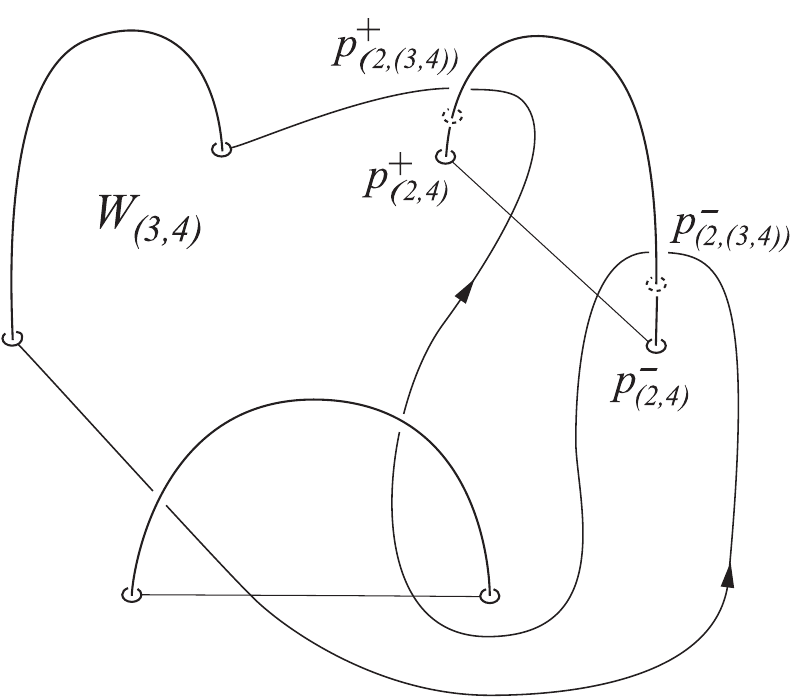}}
         \caption{}
         \label{IHX-3A-fig}
\end{figure}
We have constructed an order 1 Whitney tower $\cW^0$ which is {\em
clean}, meaning that $\cW^0$ has no unpaired intersection points
and hence is in fact a Whitney tower of order $n$ for all $n$. As
illustrated in Figure~\ref{IHX-1+3-fig}, the three order~1 Whitney
disks of $\cW^0$ all lie in the present slice of local
coordinates. In the following construction, these three Whitney
disks will be modified to create the three terms in the IHX
relation. The modified $W_{(3,4)}$ will remain entirely in the
present, while most of $W_{(2,4)}$ will be perturbed slightly into
the future, and most of $W_{(4,1)}$ will be perturbed slightly
into the past. These perturbations are essential for keeping the
Whitney disks disjoint!

Continuing with the construction, change $W_{(3,4)}$ by isotoping
its boundary $\partial W_{(3,4)}$ along $A_4$ and across
$p^+_{(2,4)}$ and $p^-_{(2,4)}$ as indicated in
Figure~\ref{IHX-3A-fig} and extending this isotopy to a collar of
$\partial W_{(3,4)}$. Note that a cancelling pair of order~1
intersection points $p^{\pm}_{(2,(3,4))}$ has been created between
$A_2$ and the interior of the ``new'' $W_{(3,4)}$ (still denoted
by $W_{(3,4)}$). The pair $p^{\pm}_{(2,(3,4))}$ is indicated in
Figure~\ref{IHX-3A-fig} by the small dashed circles near
$p^{\pm}_{(2,4)}$ and, since the orientation of $A_4$ is the
standard counter-clockwise orientation of the plane, the sign of
$p^+_{(2,(3,4))}$ (resp. $p^-_{(2,(3,4))}$) agrees with the sign
of $p^+_{(2,4)}$ (resp. $p^-_{(2,4)}$). By perturbing most of
$W_{(2,4)}$ into the future, we may assume that
$p^{\pm}_{(2,(3,4))}$ lie near, but not on, $\partial W_{(2,4)}$.
Specifically, the only part of $W_{(2,4)}$ that we do not push
into the future is a small collar of the arc of $\partial
W_{(2,4)}$ which lies on $A_4$. For now, $W_{(3,4)}$ has
intersections with the other first order Whitney disks in and near
its boundary on $A_4$, but these will be removed later in the
construction.
\begin{figure}
\subfigure[]{\includegraphics[scale=.50]{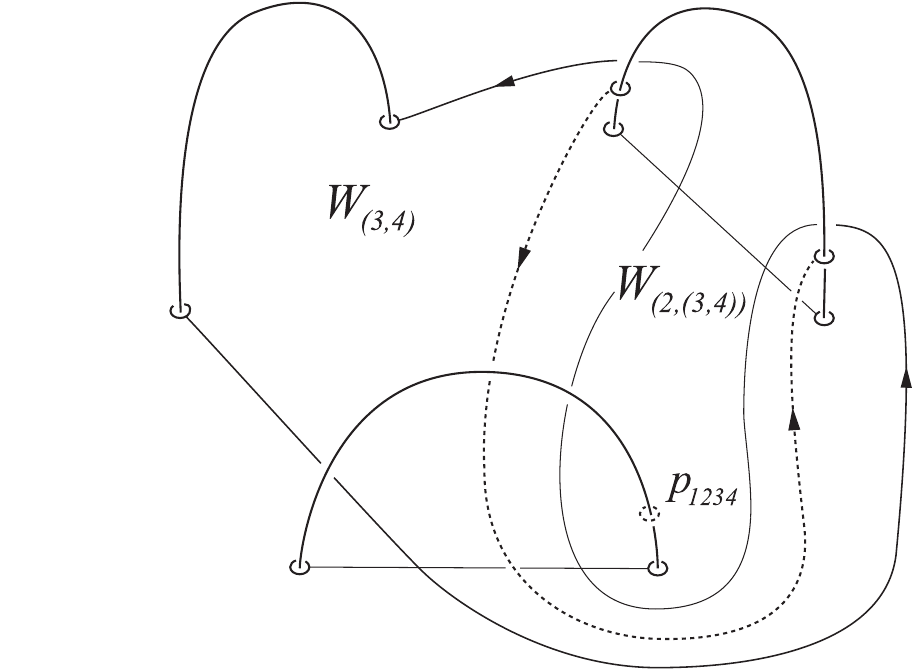}}\hfill
\subfigure[]{\includegraphics[scale=.50]{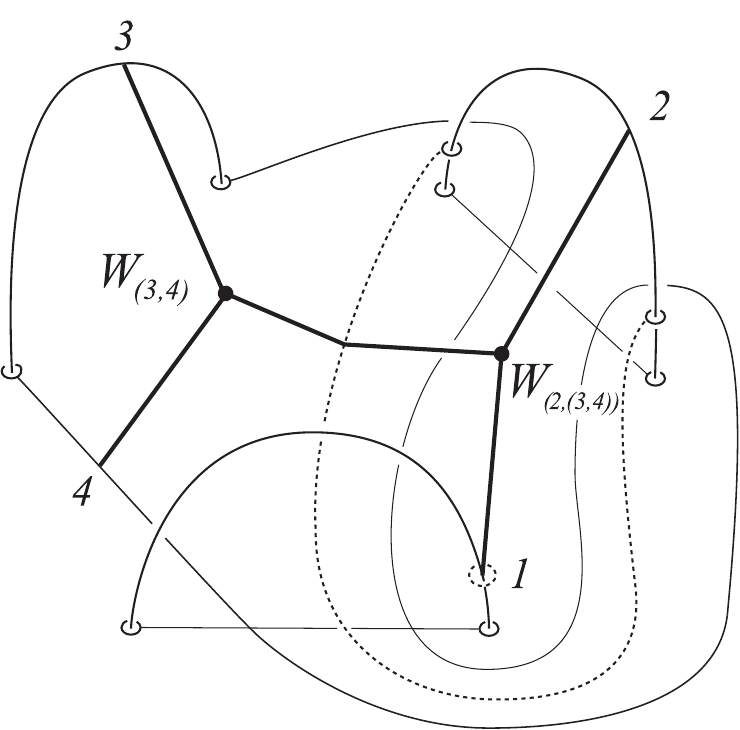}} \caption{The
construction of the $I$-tree. Both sides (a) and (b) of this
figure show the same present slice of local
coordinates.}\label{IHX-4-5-fig}
\end{figure}

A Whitney disk $W_{(2,(3,4))}$ (of order 2) for the cancelling
pair $p^{\pm}_{(2,(3,4))}$ can be constructed by altering a
parallel copy of $W_{(2,4)}$ in a collar of its boundary as
indicated in Figure~\ref{IHX-4-5-fig}a. Note that $W_{(2,(3,4))}$
sits entirely in the present. The part of the boundary of
$W_{(2,(3,4))}$ that lies on $W_{(3,4)}$ is indicated by a dashed
line in Figure~\ref{IHX-4-5-fig}a. The other arc of $\partial
W_{(2,(3,4))}$ runs along $A_2$ where there used to be an arc of
$\partial W_{(2,4)}$ before most of $W_{(2,4)}$ was pushed into
the future.

Take the orientation of $W_{(2,(3,4))}$ that corresponds to its
bracket sub-script via the above convention, i.e., that induced by
orienting $\partial W_{(2,(3,4))}$ from $p^-_{(2,(3,4))}$ to
$p^+_{(2,(3,4))}$ along $A_2$ and from $p^+_{(2,(3,4))}$ to
$p^-_{(2,(3,4))}$ along $W_{(3,4)}$ together with a second inward
pointing vector.

Note that $W_{(2,(3,4))}$ has a single {\em positive} intersection
point $p_{1234}$ (of order~2) with $A_1$ (in the present). By
pushing most of $W_{(4,1)}$ into the past, we can arrange that
$W_{(2,(3,4))}$ (which sits entirely in the present) is disjoint
from $W_{(4,1)}$. Specifically, the only part of $W_{(4,1)}$ that
is not pushed into the past is a small collar on the arc of
$\partial W_{(4,1)}$ which lies in $A_4$. To the point $p_{1234}$
we associate the positively signed labeled $I$-tree (of order~2)
as illustrated in Figure~\ref{IHX-4-5-fig}b. This $I$-tree,
$t({p_{1234}})$, is embedded in the construction with the
trivalent vertices lying in the interiors of the Whitney disks,
$W_{(2,4)}$ and $W_{(2,(3,4))}$, and each $i$-labeled univalent
vertex lying on $A_i$. Each trivalent vertex of $t({p_{1234}})$
inherits a cyclic orientation from the ordering of the components
in the bracket associated to the corresponding oriented Whitney
disk. Note that the pair of edges which pass from a trivalent
vertex down into the lower order surfaces paired by a Whitney disk
determine a ``corner'' of the Whitney disk which does not contain
the other edge of the trivalent vertex. If this corner contains
the {\em positive} intersection point paired by the Whitney disk,
then the vertex orientation and the Whitney disk orientation agree
(\cite{ST1}). Our figures are all drawn to satisfy this
convention.

\begin{figure}
\subfigure[]{\includegraphics[scale=.50]{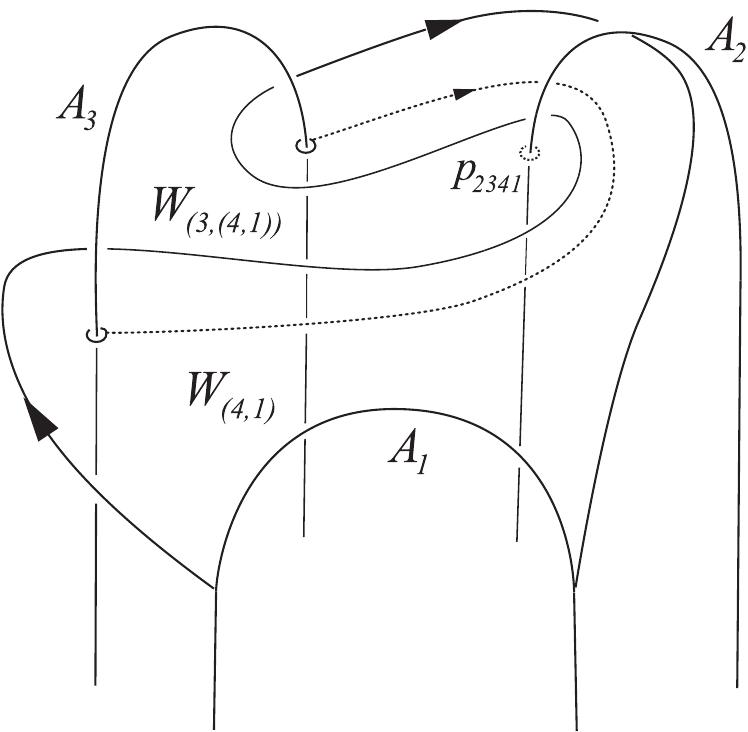}}\hfill
\subfigure[]{\includegraphics[scale=.50]{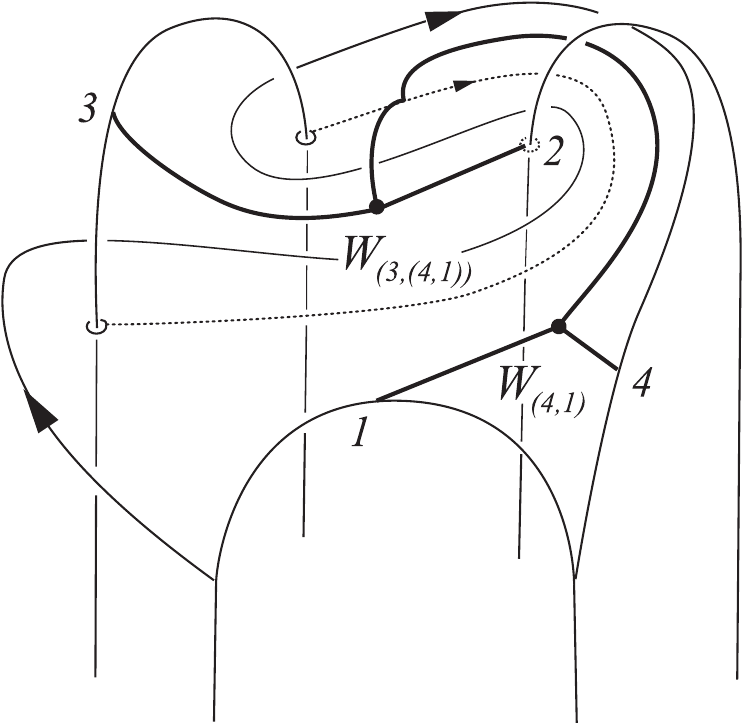}}
\caption{The construction of the $H$-tree. Both (a) and (b) of
this figure show the same slice of local coordinates, just in the
past of Figure~\ref{IHX-4-5-fig}.}\label{fig:the-H-tree-AandB}
\end{figure}

We have described how to construct (from the original $W_{(3,4)}$
of $\cW^0$) Whitney disks $W_{(2,(3,4))}$ and $W_{(3,4)}$ (both
lying entirely in the present) such that $W_{(2,(3,4))}$ pairs
$A_2\cap W_{(3,4)}$ and such that $A_1\cap W_{(2,(3,4))}$ consists
of a single point $p_{1234}$ whose associated tree is the $I$ term
in the IHX relation. In fact, a parallel version of this
construction can be carried out simultaneously on all of the
original Whitney disks in $\cW^0$ yielding additional order~2
intersection points $p_{2341}\in A_2\cap W_{(3,(4,1))}$ (with
negative sign and associated labeled trivalent tree $H$) and
$p_{3124}\in A_3\cap W_{(1,(2,4))}$ (with positive sign and
associated labeled trivalent tree $X$). Here $W_{(3,(4,1))}$ pairs
$A_3\cap W_{(4,1)}$ and $W_{(1,(2,4))}$ pairs $A_1\cap W_{(2,4)}$
and it can be arranged that all the Whitney disks have pairwise
disjointly embedded interiors and pairwise disjointly embedded
boundaries: To see this, first observe that the boundaries of the
first order Whitney disks $W_{(3,4)}$, $W_{(4,1)}$ and $W_{(2,4)}$
can be disjointly embedded in the present, as pictured in
Figure~\ref{IHX-7B-fig}, which shows how collars on the parts of
the Whitney disk boundaries that lie on $A_4$ can be
simultaneously changed in the same way that we previously changed
$W_{(3,4)}$. Recall that in the above construction, the part of
$W_{(4,1)}$ that was pushed into the past was exactly the
complement of a collar on the boundary arc of $\partial W_{(4,1)}$
which lies on $A_4$. Thus, (a collar on) the boundary arc of
$\partial W_{(4,1)}$ which lies on $A_4$ as pictured in
Figure~\ref{IHX-7B-fig} can be extended (without creating any new
intersections) to connect to the rest of $W_{(4,1)}$ which has
been perturbed into the past, and the $-H$ term can be created by
a parallel construction to the construction of the $I$ term, as
illustrated in Figure~\ref{fig:the-H-tree-AandB} which shows the
relevant past slice of local coordinates. Specifically, the second
order Whitney disk $W_{(3,(4,1))}$ sits entirely in the past, and
is made from a parallel copy of $W_{(3,4)}$ by pushing a collar to
create the intersection $p_{2341}$ with $A_2$. Note that since
$A_4$ sits entirely in the present, it does not appear in
Figure~\ref{fig:the-H-tree-AandB} which shows exclusively the
past. The signs of all intersection points can be determined from
the signs of the original intersections in
Figure~\ref{IHX-1+3-fig} using our orientation conventions: The
vertex orientations of the embedded $H$-tree in
Figure~\ref{fig:the-H-tree-AandB}(b) agree with the orientations
of the Whitney disks, and the sign of the intersection point
$p_{2341}$ is $-1$, as desired.

\begin{figure}
\subfigure[]{\includegraphics[scale=.50]{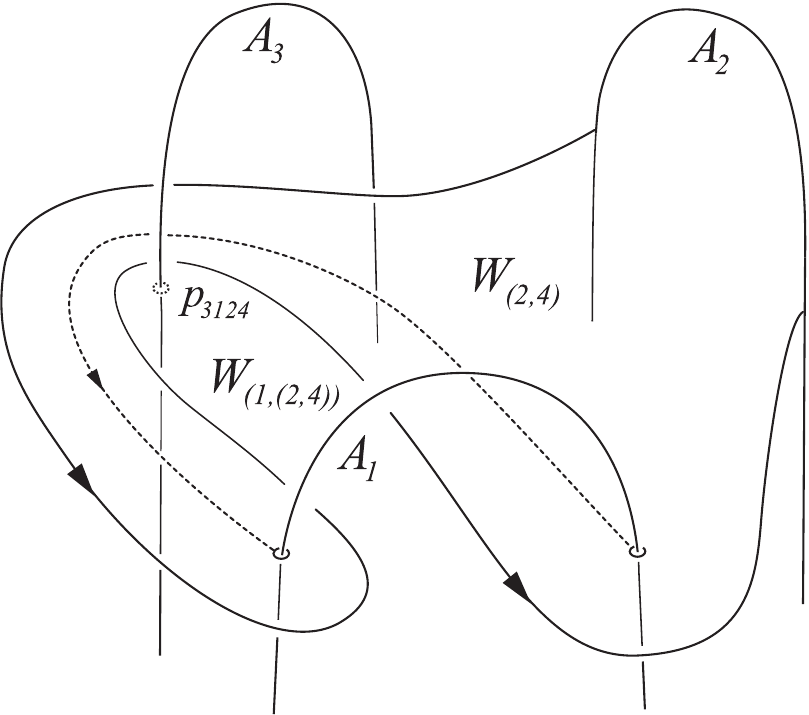}}\hfill
\subfigure[]{\includegraphics[scale=.50]{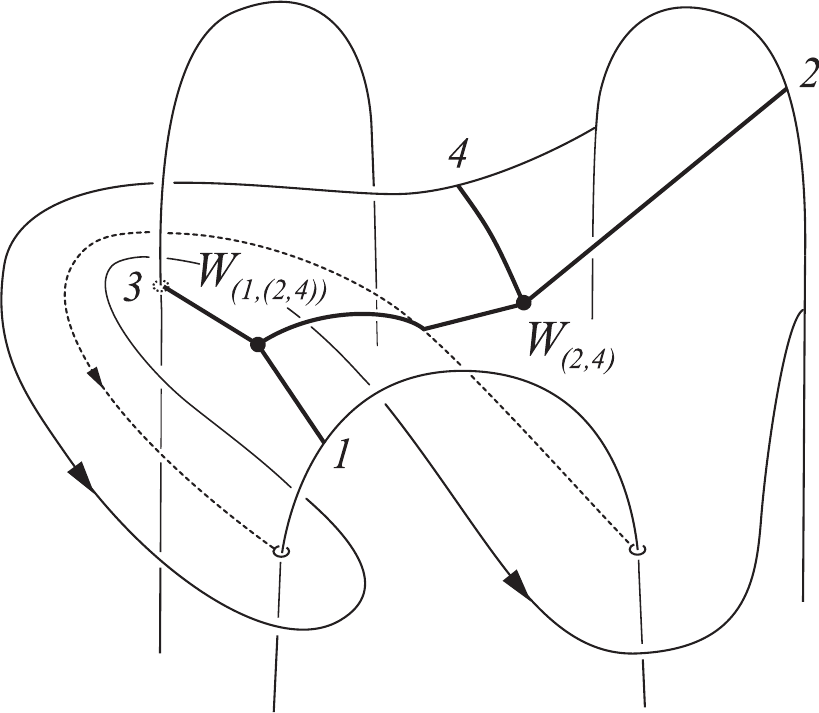}}
\caption{The construction of the $X$-tree. Both (a) and (b) of
this figure show the same slice of local coordinates, just in the
future of Figure~\ref{IHX-4-5-fig}.}\label{fig:the-X-tree-AandB}
\end{figure}

\begin{figure}[ht!]
         \centerline{\includegraphics[scale=.5]{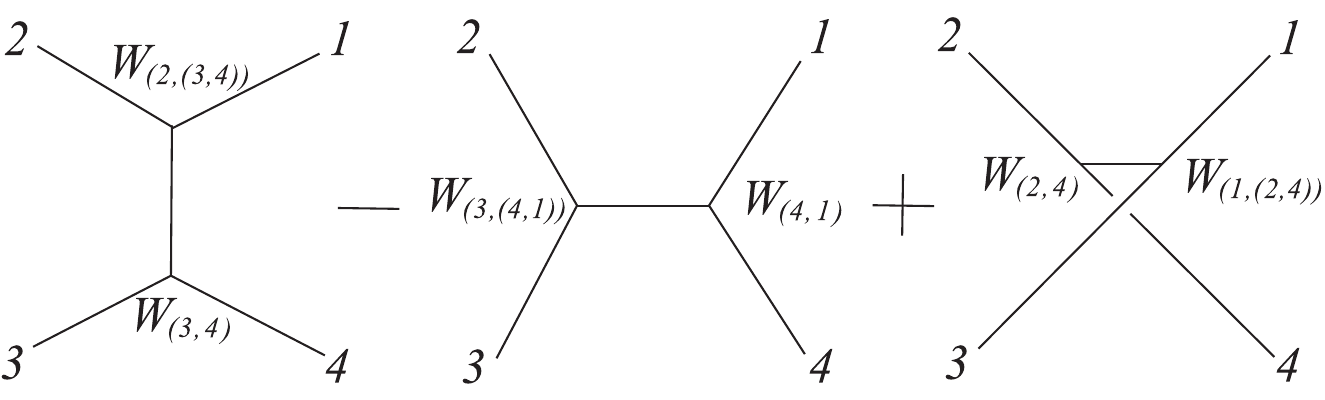}}
         \caption{The correspondence between the trivalent vertices in the
         IHX relation and the (oriented) Whitney disks in the construction.
         (The trivalent orientations are all counter-clockwise.)}
         \label{IHX-relation-1234-B-fig}
\end{figure}

The $X$-tree term is created similarly by extending a collar of
the boundary arc of $\partial W_{(2,4)}$ as pictured in
Figure~\ref{IHX-7B-fig} into the future and performing a parallel
construction as illustrated in Figure~\ref{fig:the-X-tree-AandB}.
The resulting order 2 Whitney tower $\cW$ has exactly three
order~2 intersection points with $\widetilde{\tau}_2(\cW)=(+I)
\amalg (-H) \amalg (+X)$. The correspondence between the Whitney
disks in this construction and the trivalent vertices in the IHX
relation is indicated in Figure~\ref{IHX-relation-1234-B-fig}.

The proof of the Main Theorem~\ref{thm:4dihx} is now complete, but
before moving on to connections with the 3-dimensional Jacobi
relations we note here a lemma which can now be appreciated by the
reader who has carefully kept track of the orientations in the
above constructions.

\begin{lem}\label{lem:W-tower-AS}
For a fixed order $n$ Whitney tower $\cW$, the geometric
intersection tree
$\widetilde{\tau}_n(\cW)\in\widetilde{\cB}^t_{n+1}(\ell)$ only
depends on the orientations of the order zero surfaces.
\end{lem}
\begin{proof}
Recall that $\widetilde{\tau}_n(\cW)$ is a disjoint union of
signed vertex oriented trees associated to the order $n$
intersection points in $\cW$, and the AS relations change the sign
of a tree whenever a vertex orientation is changed. Each tree
$t(p)$ is most easily defined as a subset of $\cW$ which
bifurcates down through the Whitney disks, with each trivalent
vertex of $t(p)$ lying in a Whitney disk. Each trivalent vertex
has two \emph{descending} edges which pass into the lower order
sheets paired by the Whitney disk, and one \emph{ascending} edge
which either passes through the intersection point $p$ or passes
into a higher order Whitney disk. Assuming fixed orientations on
all the surfaces in $\cW$ (including Whitney disks),  our
orientation convention for $t(p)$ can be summarized as follows:
The descending edges of a trivalent vertex determine a corner of
the corresponding Whitney disk which does not contain the
ascending edge. If this corner encloses the positive intersection
point (of the intersections paired by the Whitney disk), then the
vertex-orientation is the same as that induced by the orientation
of the Whitney disk. If this corner encloses the negative
intersection point, then the vertex-orientation is the opposite of
the orientation induced by the Whitney disk.

We remark here that in practice the geometric intersection tree
$\widetilde{\tau}_n(\cW)$ usually sits as an \emph{embedded}
subset of $\cW$, as can be arranged easily by ``splitting'' the
Whitney tower (\cite{S,ST1}). However in general
$\widetilde{\tau}_n(\cW)$ will not be embedded if any Whitney
disks contain self-intersections and/or multiple (pairs of)
intersections.

To check that each signed tree $\e _p \cdot t(p)$ in
$\widetilde{\tau}_n(\cW)$ only depends, modulo antisymmetry, on
the orientations of the underlying order zero surfaces it is
enough to consider the effect of changing any single Whitney disk
orientation. There are two cases to consider:

First consider a Whitney disk $W$ containing a trivalent vertex
$v$ of a signed tree $\e _p \cdot t(p)$, where the ascending edge
of $v$ passes into a higher order Whitney disk $W'$ containing an
adjacent trivalent vertex $v'$ of $t(p)$. Changing the orientation
of $W$ changes the vertex-orientation of $v$, and also changes the
vertex-orientation at $v'$ because the signs of the intersection
points (in $W$) which are paired by $W'$ are reversed. Thus, the
signed tree $\e _p \cdot t(p)$ does not change.

Now consider a Whitney disk $W$ containing a trivalent vertex $v$
of a signed tree $\e _p \cdot t(p)$, where the ascending edge of
$v$ passes through the intersection point $p$. In this case,
changing the orientation of $W$ changes the vertex-orientation of
$v$ and changes the sign of the intersection point $p$, provided
$p$ is not a self-intersection point of $W$, so that $\e _p \cdot
t(p)$ is changed exactly by an antisymmetry relation at $v$. In
the case that $p$ is a self-intersection point of $W$, then
changing the orientation of $W$ changes both trivalent vertices
adjacent to $p$, namely $v$ and another trivalent vertex of $t(p)$
which also sits in $W$.
\end{proof}

\section{Connecting 4- and 3-dimensional Jacobi relations}\label{sec:3and4}
In this section we explain in detail the commutative diagram
(\ref{diagram}) in \ref{sec:down-to-3} of the introduction. But
first we need to introduce some background material.

\subsection{Gropes and their associated trees}\label{sec:gropes}

For technical simplicity, we will use only genus one gropes, which
are sufficient for our purposes. We will not specify the genus one
assumption in the body of this paper but we note that there is a
{\em grope refinement} procedure \cite[2.3]{CT1} that allows one
to replace an arbitrary grope by a genus one grope. In fact, we
allow here the bottom stage surface to have arbitrary genus,
that's why we don't need {\em sequences of} genus one gropes as in
\cite{CT1}. In the appendix we deal with general gropes, and the
reader is referred to \cite{CT1} for their definition.

\begin{defi}\label{def:genus-one-gropes}
A (genus one) {\em  grope} $g$ is constructed by the following
method:
\begin{itemize}
\item Start with a compact oriented connected surface of any
genus, the {\em bottom stage} of $g$, and choose a symplectic
basis of circles on this bottom stage surface.
\item Attach punctured tori to any number of the basis circles and
choose hyperbolic pairs of circles on each attached torus.
\end{itemize}
Iterating the second step a finite number of times yields the
grope $g$. The attached tori are the {\em higher stages} of $g$.
The basis circles in all stages of $g$ that do not have a torus
attached to them are called the {\em tips} of $g$. Attaching
2--disks along all the tips of $g$ yields a {\em capped} grope (of
genus one), denoted $g^c$. In the case of an (uncapped) grope, it
is often convenient to attach an annulus along one of its boundary
components to each tip. These annuli are called \emph{pushing
annuli}, and every tame embedding of a grope in a 3--manifold can
be extended to include the pushing annuli. \hfill$\Box$
\end{defi}

Let $g^c$ be a capped  grope. We define a rooted trivalent tree
$t(g^c)$ as follows:
\begin{defi}\label{grope tree}
Assume first that the bottom stage of $g^c$ is a genus one surface
with boundary.  Then define $t(g^c)$ to be the rooted trivalent
tree which is dual to the 2--complex $g^c$; specifically, $t(g^c)$
sits as an embedded subset of $g^c$ in the following way: The root
univalent vertex of $t(g^c)$ is a point in the boundary of the
bottom stage of $g^c$, each of the other univalent vertices are
points in the interior of a cap of $g^c$, each higher stage of
$g^c$ contains a single trivalent vertex of $t(g^c)$, and each
edge of $t(g^c)$ is a sheet-changing path between vertices in
adjacent stages or caps (here ``adjacent'' means ``intersecting in
a circle''), see Figure~\ref{posquad12}b.

In the case where the bottom stage of $g^c$ has genus $>1$, then
$t(g^c)$ is defined by cutting the bottom stage into genus one
pieces and taking the disjoint union of the trees just described.
In the case of genus zero, $t(g^c)$ is the empty tree.\hfill$\Box$
\end{defi}

We can now define the relevant complexity of a grope.
\begin{defi}\label{grope class}
The {\em class} of $g^c$ is the minimum of the Vassiliev degrees
of the connected trees in $t(g^c)$. The underlying uncapped grope
$g$ (the {\em body} of $g^c$) inherits the same tree,
$t(g)=t(g^c)$, and the same notion of class. If the grope consists
of a surface of genus zero, we regard it as a grope of class $n$
for all $n$. The non-root univalent vertices of $t(g)$ are called
{\em tips} and each tip of $t(g)$ corresponds to a tip of
$g$.\hfill$\Box$
\end{defi}
We will assume throughout the paper that all surface stages in our
gropes contribute to the class of the grope, i.e. we ignore
surface stages that can be deleted without changing the class.

\subsection{Grope cobordism} \label{sec:grope cobordism}

\begin{defi}\label{grope cobordism}
A {\em class n grope cobordism} between $\ell$-component string
links $\sigma$ and $\sigma'$ is defined as follows. For each
$1\leq i\leq \ell$, let $\sigma_i$, respectively $\sigma_i'$, be
the $i$-th string link component of $\sigma$, respectively
$\sigma'$. Suppose that, for each $i$, there is an embedding of a
class $n$ grope $g_i$ into the 3-ball whose (oriented) boundary is
decomposed into two arcs representing the (oriented) isotopy
classes of $\sigma_i$ and $-\sigma_i'$. This collection of gropes
is called a \emph{grope cobordism} $G$ from $\sigma$ to $\sigma'$
if the gropes $g_i$ are embedded disjointly. We sometimes also say
that $G$ is a grope cobordism {\em of} $\sigma$ and note the
asymmetry coming from the above orientation convention.

If all the tips of  each $g_i$ bound embedded caps whose interiors
are disjoint from each other and disjoint from all but the bottom
stages of the $g_i$, then $G$ together with these caps forms a
(class $n$) {\em capped} grope cobordism $G^c$ from $\sigma$ to
$\sigma '$ (or {\em of} $\sigma$).
 \hfill$\Box$
\end{defi}
Note that this definition does not specify the relative embedding
of $\sigma$ and $\sigma'$.
\begin{rem}
The above definition is a generalization of the one given in
\cite{CT1} for knots. By considering disjointly embedded gropes in
3-space, each with two boundary circles, one also gets a notion of
grope cobordism of links. The arguments of \cite{CT1} adapt to
show that grope cobordism (of links or string links) is an
equivalence relation. \hfill$\Box$
\end{rem}

Let $G^c$ be a  capped grope cobordism from $\sigma$ to $\sigma'$.
It turns out that one can assume that the intersections of the
caps with the bottom stages are arcs from $\sigma$ to
$\sigma^\prime$. This can be accomplished by finger moves of the
caps across the boundary of the bottom stages. Also, by applying
Krushkal's splitting technique (as adapted to 3--dimensions in
\cite{CT1}) it can be assumed that each cap contains just a single
intersection arc.

\begin{defi} The following notions will be used for capped grope cobordisms.
\begin{enumerate}
\item A  capped grope cobordism which has been simplified
as described in the previous paragraph will be referred to as a
\emph{simple grope cobordism.}
\item Denote by $\G^c_n(\ell)$ the set of
class $n$  simple grope cobordisms of $\ell$-component string
links.
\item Denote by $\G_n(\ell)$ the set of class $n$  grope cobordisms.
(That is,  grope cobordisms which are not required to have
caps.)\hfill$\Box$
\end{enumerate}
\end{defi}

\subsection{Claspers and gropes} \label{sec:claspers}

\begin{defi} The following definitions can be found in
\cite{H} and/or \cite{CT1}.
\begin{enumerate}
\item A \emph{clasper} is a surface embedded in the complement of a link or string link in a $3$-manifold,
formed by gluing together edges, nodes and leaves. An edge is
homeomorphic to $I\times I$, and each end $I\times\{0\}$ or
$I\times\{1\}$ is glued to a node or a leaf. A node is
homeomorphic to $D^2$ and must have three edges glued to its
boundary. A leaf is homeomorphic to $S^1\times I$ and must have a
singe edge glued to one of its boundary components.

\item A clasper is said to be \emph{capped} if all of (the cores of)
its leaves bound disjoint disks (called caps) which may hit the
link or string link, but only intersect the clasper along their
boundaries.

\item A clasper is said to be \emph{simple} if it is capped and if the caps each only hit the
link or string link in a single transverse intersection.

\item Given a clasper $C$, we can form an oriented graph by collapsing each edge to a $1$-dimensional edge,
each node to a trivalent vertex, and each leaf to a univalent
vertex. The vertex-orientation of the graph is somewhat subtle,
especially when the resulting graph is not a tree, and we refer
the reader to \cite{CT2} for details.

\item A \emph{tree} clasper is a clasper whose associated graph is a tree.

\item A tree clasper is said to be \emph{rooted} if there is at least one leaf which has a
cap that hits the link or string link in a single transverse
intersection.

\item  Given a clasper, there is a way of producing an embedded framed link,
and surgery on the clasper is defined to be surgery on this framed
link. If the clasper is rooted (which is implied by simple and
capped) then the surgery does not change the ambient manifold and
can instead be regarded as changing the link or string link.
\hfill$\Box$
\end{enumerate}
\end{defi}

\begin{defi}
The (Vassiliev) \emph{degree} of a clasper is half the total
number of vertices of the associated graph. The \emph{grope
degree} of a clasper is the (Vassiliev) degree plus the first
Betti number of the associated graph.\hfill$\Box$
\end{defi}

Claspers and gropes are closely related, as discussed in detail in
\cite{CT1}. Here are some important results, which were stated for
knots, but hold true for links and string links as well.

\begin{thm} The following statements can be proven by the techniques of \cite{CT1}.
\begin{enumerate}
\item Two links or string links in a $3$-manifold differ by a sequence of
simple clasper surgeries of Vassiliev degree $n$ if and only if
they are related by a simple grope cobordism of class $n$.
\item Two links or string links in a $3$-manifold differ by a sequence of rooted
tree clasper surgeries of Vassiliev degree $n$ if and only if they
are related by a grope cobordism of class $n$.
\item Two links or string links in a $3$-ball differ by a sequence of simple clasper
surgeries of grope degree $n$ if and only if they are related by a
grope cobordism of class $n$.
\end{enumerate}
\end{thm}

Habiro \cite{H} has shown that two knots share the same Vassiliev
invariants up to degree $n$ if and only if they differ by a
sequence of simple clasper surgeries of degree $(n+1)$. Together
with the above Theorem, this implies two knots have the same
Vassiliev invariants up to degree $n$ if and only if they cobound
a simple grope cobordism of class $(n+1)$. The corresponding
statements for string links are not known, but see
Section~\ref{sec:string}.

\subsection{Geometric intersection trees for grope
cobordisms}\label{subsec:geo-int-trees-for-grope-cobordism} Let
$G^c\in\mathbb G^c_n(\ell)$ be a class $n$ simple grope cobordism
of a string link $\sigma$, and let $g_i^c$ be a capped grope
component of $G^c$.  Each cap of $g_i^c$ contains only a single
arc of intersections, which can be with any bottom stage surface
in $g_j^c\subset G^c$. The bottom stage surface of $g_i^c$
inherits an orientation from its boundary, and we now describe how
to orient the higher stages of the grope cobordism, up to a
certain indeterminacy.

\begin{figure}[ht!]
\centering \subfigure[Positive quadrants and
orientation.]{\includegraphics[width=.4\linewidth]{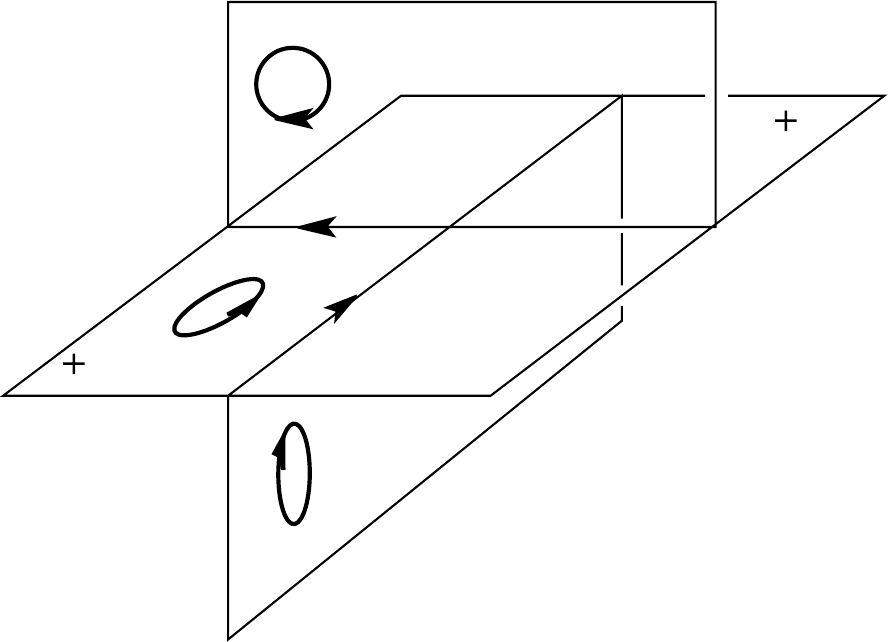}}
\subfigure[A trivalent vertex of
$t(g^c)$.]{\includegraphics[width=.4\linewidth]{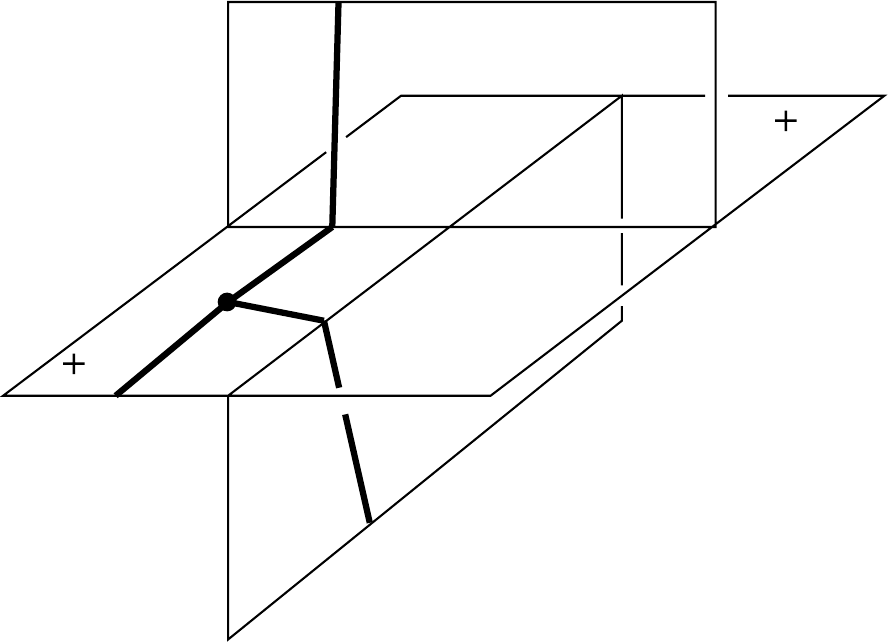}}
\caption{}\label{posquad12}
\end{figure}

 Each surface stage or cap is attached
to a previous stage along a circle, which hits the attaching
region for one other surface stage or cap in a point. Near this
point, the 2-complex is modeled by the following subset of $\R^3$:
$$\{(x,y,z):z=0\}\cup \{(x,y,z): x=0, z\geq 0\} \cup\{ (x,y,z): y=0, z\leq 0\}.$$
Distinguish two of the quadrants as \emph{positive}, namely the
quadrants where both $x,y>0$ respectively where both $x,y <0$. See
Figure~\ref{posquad12}a, where the two positive quadrants are
indicated.
Now suppose that the lower stage $(z=0)$ has an orientation and
choose one of the two positive quadrants. The orientation of the
surface induces an orientation of a small triangle in the positive
quadrant which has a vertex at the origin and two edges contained
in the axes. This then induces an orientation of the boundaries of
the two higher surface stages, and hence induces an orientation of
the higher surface stages.
 If we use the other positive quadrant instead, this
has the effect of flipping the orientation of both higher surface
stages, and this is the indeterminacy that we allow.

The above orientations of the surface stages in a capped grope
$g^c$ induce vertex-orientations of the trivalent vertices of
$t(g^c)$ by taking each trivalent vertex of $t(g^c)$ to lie in a
chosen positive quadrant, see Figure~\ref{posquad12}b. Here also,
the pairs of edges that cross into the next stages are required to
do so through that positive quadrant.

Recall that $t(g_i^c)$ is a disjoint union $\amalg_r t_i^r$ of
trees $t_i^r$, each of which sits as an embedded subset of
$g_i^c$, with the root of $t_i^r$ lying on the $i$-th strand of
$\sigma$ (in $\partial g_i^c$) and each tip of $t_i^r$ lying
inside a cap. The interior of each cap intersects the cobordism in
a single intersection arc which corresponds to some strand of
$\sigma$. Hence we can regard these tips as actually lying on
 a $j$-th strand of $\sigma$ at an
intersection point between a cap of $g_i^c$ and that $j$-th strand
(see left hand side of
Figure~\ref{oriented-grope-surgery-with-trees3R-fig}). Associate
to each tip of $t_i^r$ the sign of the corresponding intersection
point (between the cap and the $j$-th strand) and denote by
$\epsilon_i^r\in\{+,-\}$ the product of these signs.

The vertices of $t(g_i^c)$ can be oriented by regarding the tree
as a subset of $g_i^c$ where the two edges emanating from a
trivalent vertex must pass to the higher stages in a positive
quadrant, as depicted in Figure~\ref{posquad12}.

Recall $\widetilde{\cA}^t_n(\ell)$ from Definition~\ref{def:A}.
\begin{defi}\label{def:cappedtau}
Let $G^c$ be a capped class $n$ grope cobordism of $\ell$-string
links. The {\em geometric intersection tree}
$\widetilde{\tau}^c_n(G^c)\in\widetilde{\cA}^t_n(\ell)$  is
defined to be the disjoint union
$\amalg_i\amalg_r\epsilon_i^r\cdot t_i^r$ of all the
vertex-oriented signed trees associated to all the $g_i^c$. Note
that each tree should avoid the intersections between caps and the
bottom stage, and this forces the roots to attach to the strands
of $\sigma$ in a specific ordering.
\end{defi}

\begin{lem} \label{lem:AS for gropes}
The {\em geometric intersection tree} $\widetilde{\tau}^c_n(G^c)$
is well-defined.
\end{lem}
\begin{proof}
The issue is whether the choice of positive quadrants can affect
$\widetilde{\tau}^c_n(G^c)$. Choosing a different positive
quadrant does not change the cyclic order of the corresponding
vertex, but it does change the orientations of all of the higher
stages, including the caps. This has the effect of switching the
cyclic orders at each of the vertices representing these higher
stages, as well as switching the sign of all of the tips
representing these caps. In other words a sign is introduced for
every vertex (both $1$- and $3$-valent) lying above the vertex we
started with. A simple induction shows that there must be an even
number of these. Hence, we arrive at the same signed tree, modulo
AS relations.
\end{proof}

\begin{defi}\label{def:Psi}
\begin{enumerate}
\item Let $\rho\colon \widetilde{\cA}^t_n(\ell)\to \widehat{\cA}_n^t(\ell)$ be the natural map
sending the monoid operation to the group addition.
\item Let $\widehat\tau^c_n\colon\G^c_n(\ell)\ra
\widehat{\cA}_n^t(\ell)$ be defined by
$\widehat\tau^c_n(G^c)=\rho(\widetilde{\tau}^c_n(G^c))$.
\end{enumerate}
\hfill$\Box$
\end{defi}

\begin{rem}
If one translates a simple grope into a union of simple tree
claspers, the map $\widehat\tau^c_n$ can be regarded as the map
which sums over the set of claspers, collapsing each clasper to
its underlying tree, with univalent vertices attaching to the
$\ell$ strands according to where the caps of the clasper meet the
string link. This was the point of view taken in \cite{CT2}.
$\hfill\Box$
\end{rem}

\subsection{From grope cobordism to Whitney concordance}\label{sec:3to4}

\begin{defi}
A {\em singular concordance} between string links $\sigma$ and
$\sigma'$ is a collection of properly immersed 2--disks $D_i$ in
the product $B^3\times I$ of the 3--ball with the unit interval
$I=[0,1]$, with $\partial D_i$ equal to the union of the $i$-th
strands $\sigma_i \subset B^3\times \{0\}$ and $\sigma_i'\subset
B^3\times \{1\}$ together with their end points crossed with $I$.
For instance, any generic homotopy between $\sigma$ and $\sigma'$
defines such a singular concordance. A singular concordance {\em
of} $\sigma$ induces the orientation of $\sigma$.

An (order $n$) Whitney tower whose bottom stages form a singular
concordance is called an (order $n$) {\em Whitney concordance}.
Denote by $\W_n(\ell)$ the set of order $n$ Whitney concordances
of $\ell$-component string links. \hfill$\Box$
\end{defi}

Let $G^c$ be a simple grope cobordism (from $\sigma$ to $\sigma'$)
in $\G^c_n(\ell)$. Think of $G^c$ as sitting in the middle slice
$B^3\times\{1/2\}$ of $B^3\times I$.  Extending $\sigma\subset
G^c$ to $B^3\times \{0\}$, via the product with $[0,1/2]$, and
extending $\sigma'\subset G^c$ to $B^3\times \{1\}$, via the
product with $[1/2,1]$, yields a collection of class $n$ capped
gropes properly embedded in $B^4=B^3\times I$, i.e. a {\em grope
concordance}, from $\sigma$ to $\sigma'$. After perturbing the
interiors of the caps slightly, we may assume that all caps are
still disjointly embedded and that a cap which intersected the
$j$-th string link component in the grope cobordism now has a
single transverse intersection point with the interior of a bottom
stage of the $j$-th grope in the grope concordance. By fixing the
appropriate orientation conventions, this construction preserves
the signs of these intersection points.

\begin{figure}[ht!]
         \centering
         \includegraphics[scale=.6]{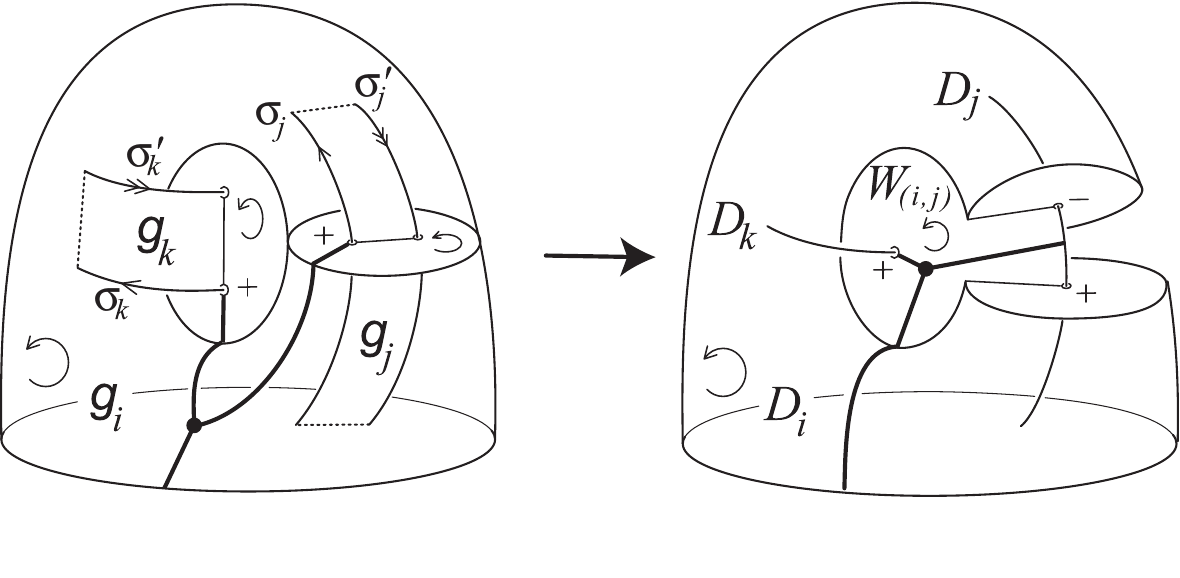}
         \caption{Left: A top stage of a capped grope cobordism. Right:
         The corresponding part of a Whitney concordance after pushing into 4--space and
         surgering a cap.}
         \label{oriented-grope-surgery-with-trees3R-fig}
\end{figure}

Consider the effect of the construction on the (degree~n) trees
$t(g_i^c)$ which were embedded in the original $G^c$ and are now
sitting in the class $n$ capped gropes in the 4--ball: Any root
vertex that was lying on an $i$-th string link strand is now in
the interior of the $i$-th bottom stage, and any tip that
corresponded to an intersection between a cap and a $j$-th strand
now corresponds to an intersection between a cap and a $j$-th
bottom stage. These are exactly the labeled  trees associated to
gropes in 4--manifolds as described in \cite{S}, and Theorem~6 of
\cite{S} describes how to surger such gropes to an order $(n-1)$
Whitney concordance $\cW$ while preserving trees, meaning that the
labeled trees associated to the gropes become the order~$(n-1)$
geometric intersection tree $\widetilde{\tau}_{n-1}(\cW)$.
Although signs and orientations are not discussed in \cite{S}, the
notation there is compatible with the sign conventions of this
paper and a basic case of the compatibility is illustrated in
Figure~\ref{oriented-grope-surgery-with-trees3R-fig} which shows a
``push and surger'' step in the modification of a 3-dimensional
grope cobordism to a Whitney concordance applied to a top stage.
The modification in general involves ``hybrid'' grope-towers but
reduces essentially to this case as explained in \cite{S}.

\begin{defi}\label{pin} The commutative diagram (\ref{diagram}) in the introduction is explained as follows.
\begin{enumerate}
\item Let $\overline{\W}_{(n-1)}$ be the set of order $(n-1)$ Whitney concordances
 modulo the relation that two Whitney towers with the same geometric intersection tree are the same.
\item The above constructions define the map $\pin: \G^c_n(\ell)  \ra
\overline{\W}_{(n-1)}(\ell) $ which pushes a grope into $B^3\times
I$ and surgers it into a Whitney tower. It is used in our main
diagram~(\ref{diagram}) in the introduction. The Whitney tower
produced from a grope is not unique, as it depends on the choice
of caps one uses to surger, which is why we need to divide
$\W_{(n-1)}(\ell)$ by an appropriate equivalence
relation.\hfill$\Box$
\end{enumerate}
\end{defi}

\begin{rem}\label{maps pq}
The only information contained in the original geometric
intersection tree $\widetilde{\tau}^c_n(G^c)$ that is {\em lost}
by the map (induced by) $\pin$ is the ordering in which the
univalent vertices of the trees in  $\widetilde{\tau}^c_n(G^c)$
were attached to the string link components. Thus, pushing a class
$n$ grope cobordism into 4--dimensions, surgering to an order
$(n-1)$ Whitney concordance and applying the map
$\widetilde{\tau}_{(n-1)}$ is the same as the composition of the
map $\widetilde{\tau}^c_n$ with the homomorphism
\[
\poff:\widehat{\cA}^t_n(\ell)\ra\widehat{\cB}^t_n(\ell)
\]
 that pulls the trees off the string link components
and labels the univalent vertices accordingly.

Notice that the map $\poff$ is very different from the rational
PBW-type isomorphism $\sigma\colon \cA\otimes\mathbb
Q\to\cB\otimes\mathbb Q$ as defined in \cite{BN}.\hfill$\Box$
\end{rem}

\section{Jacobi Identities in Dimension 3}\label{sec:3-dim}
As a consequence of our work so far,  IHX relations appear in
$\widetilde{\cB}^t_n$, and hence $\widehat{\cB}^t_n$, as the image
under $\widetilde{\tau}_{(n-1)}$ (respectively
$\widehat{\tau}_{(n-1)}$) of Whitney concordances from any
string-link to itself (e.g., tube the 2--spheres in
Theorem~\ref{thm:4dihx} into a product concordance). In Sections
\ref{sec:ihxstringlinks} and \ref{sec:genihx} we show that this
phenomenon pulls back to the 3--dimensional world: There are
capped grope cobordisms from any string link to itself whose
images under ${\widetilde{\tau}^c_n}$ (and $\widehat{\tau}^c_n$)
give all IHX relations. We will also realize all IHX relations in
a group generated by unitrivalent {\em graphs} by defining a more
general map $\widehat\tau^g_n$ on {\em un}-capped  class $n$ grope
cobordisms.  In the appendix, we will show how to interpret this
map for grope cobordisms where genus is allowed at all stages.

\subsection{The IHX relation for string links}\label{sec:ihxstringlinks}
The geometric IHX construction for string links contained in
Theorem~\ref{thm:3dihx}  will play a key role in all subsequent
IHX constructions. At the heart of the proof of
Theorem~\ref{thm:3dihx} is a 3-dimensional interpretation of
Figure~\ref{IHX-7B-fig} which leads to the following construction
of a capped grope cobordism which is (slightly) singular -- these
singularities will be removed in subsequent constructions.

\begin{construct}\label{construct:sing-grope}
Consider a trivial three-component string link in the 3--ball. We
will construct a singular capped grope $\bar{g}^c$ of class three
with an unknotted boundary component on the surface of the ball.
Its bottom stage is of genus three and embedded. The second stage
surfaces of $\bar{g}^c$ are of genus one and are each embedded.
The interiors of the second stage surfaces intersect each other
but are disjoint from the bottom stage of $\bar{g}^c$. Only the
caps of $\bar{g}^c$ intersect the three trivial string link
strands. Denote by $\bar{G}^c$ the union of $\bar{g}^c$ together
with trivial cobordisms of the strands of the string link
(embedded 2--disks traced out by perturbations of the interiors of
the strands). Then the key property of $\bar{g}^c$ is that
$\widetilde\tau^c_3(\bar{G}^c)\in\widetilde{\cA}^t_3(4)$ is equal
to the three terms of the IHX relation in
Figure~\ref{IHX-relation}. Here the strands of the trivial string
link are labeled by $1$, $2$, $3$ and $\bar{g}^c$ is interpreted
as a {\em null} bordism of its unknotted boundary which is labeled
$4$. Note that $\widetilde\tau^c_3$ still makes sense as a
disjoint union of subtrees of $\bar{g}^c$ whose tips are attached
to intersections with caps, even though $\bar{g}^c$ is singular.

To begin the construction of $\bar{g}^c$, consider
Figure~\ref{IHX-7B-fig} again. Think of it as taking place inside
a $3$-ball $B$, so that the horizontal plane has an unknotted
boundary on $\partial B$. The arcs that each puncture the plane
twice are the three strands of a trivial string link. Add tubes
around the arcs to turn the plane into a genus three surface
$\Sigma$. $\Sigma$ is the bottom stage of our singular grope
$\bar{g}^c$. We construct a symplectic basis for $\Sigma$ as
follows. Three of the curves are meridians to the tubes. To get
the other three basis curves, connect the endpoints of each of the
three pictured arcs in the plane (formerly Whitney arcs) by an
untwisted arc that travels once over a tube. (Exercise: these
three curves form a Borromean rings.) We fix surfaces bounding
these latter three basis curves in the following way. Consider
Figure~\ref{IHX-3A-fig}, where a Whitney disk $W_{(3,4)}$ is
pictured. Thinking of the figure as being in a 3-ball (rather than
a 3-dimensional slice of 4-space), the Whitney disk has two
intersections with an arc in strand 2 of the trivial string link,
and adding a tube around this arc yields a surface $s_1$ as
illustrated in Figure~\ref{fig:SingularGrope
Construction1B-and-with-tree}(a). This surface has a pair of dual
caps, whose boundaries are indicated by the dashed loops in
Figure~\ref{fig:SingularGrope Construction1B-and-with-tree}. One
of these caps intersects the upper right strand 2, and the other
intersects the bottom strand 1; these caps also have circles of
intersection (not shown in the figure) with the tubes of $\Sigma$
around these strands (but these circles of intersections will be
eliminated during later applications of this construction). The
curve dual to the attaching curve of $s_1$ is a meridian to the
strand 3 and so bounds a cap hitting strand 3 once. The tree
structure for the stage $s_1$ and its dual cap is $[[1,2],3]$, as
shown in Figure~\ref{fig:SingularGrope
Construction1B-and-with-tree}(b).
\begin{figure}
\subfigure[]{\includegraphics[scale=.50]{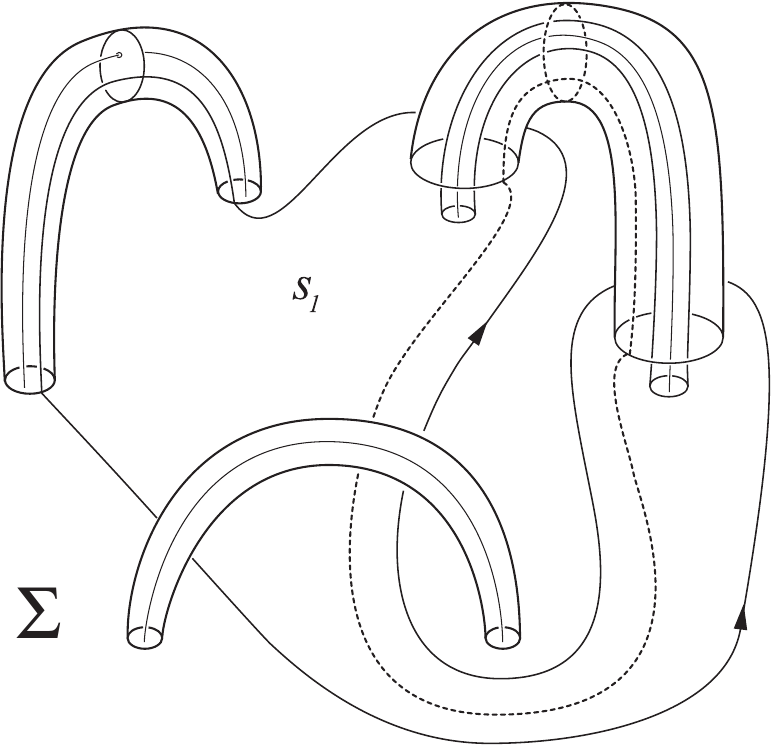}}\hfill
\subfigure[]{\includegraphics[scale=.50]{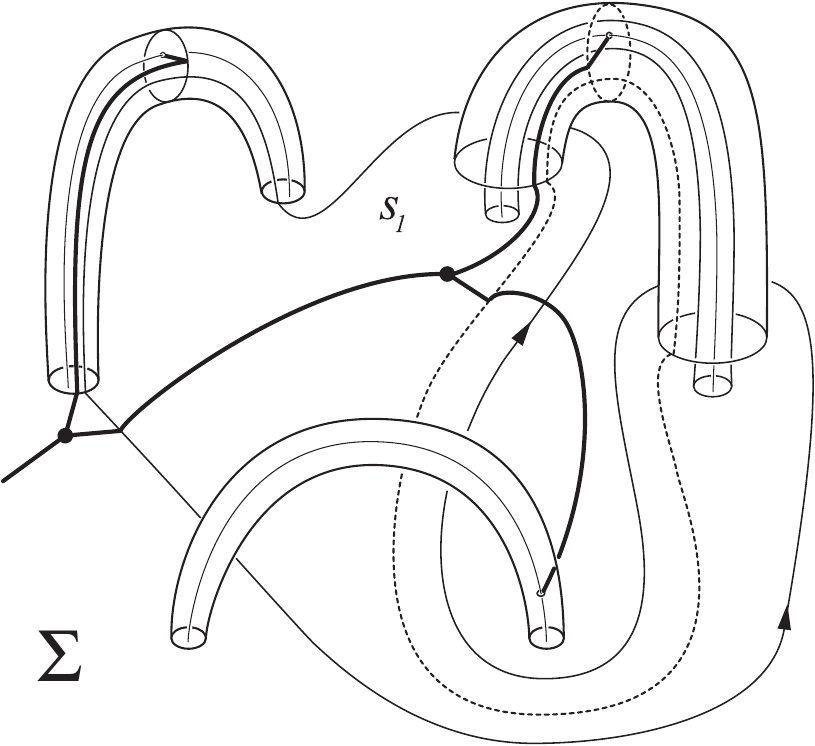}}
\caption{The construction of the capped surface $s_1$ for the
singular capped grope $\bar{g}^c$ in
construction~\ref{construct:sing-grope}.}\label{fig:SingularGrope
Construction1B-and-with-tree}
\end{figure}

Symmetrically, we can construct $s_2$ and $s_3$, with trees
$[1,[2,3]] $ and $[[3,1],2]]$, by interpreting
Figure~\ref{fig:the-H-tree-AandB} and
Figure~\ref{fig:the-X-tree-AandB} as both being in the 3-ball.
Adding these three capped surfaces $s_1,s_2,s_3$ to the surface
$\Sigma$ we get the desired singular capped grope $\bar{g}^c$
bounded by strand 4. Including strand 4 as the root, the
associated three trees give exactly the terms of the IHX relation.
With a little extra effort in analyzing the orientations, one can
verify that the signs of these three terms are correct. \hfill
$\Box$
\end{construct}

\begin{proof}[Proof of Theorem~\ref{thm:3dihx}]
First consider the case where $\ell=4$ and
$(+t_I)\amalg(-t_H)\amalg(+t_X)$ is as in
Figure~\ref{IHX-relation}. We will construct $G^c$ as a grope
cobordism of strand $4$ together with trivial cobordisms (disks)
of the other three strands. Take the 3--ball $B$ from the above
Construction~\ref{construct:sing-grope} and remove regular
neighborhoods of the three strands of the trivial string link in
$B$ to get a handlebody $M$ which contains the uncapped body
$\bar{g}$ of the singular capped grope $\bar{g}^c$. Let $m_i$ be a
meridian to the $i$-th strand on the surface of $M$. Now in the
complement of a trivial $4$-component string link, embed $M$ so
that $m_i$ is a meridian to strand $i$. Connect a parallel copy of
the fourth strand by a band to the unknot $\partial\bar{g}$ on the
boundary surface of $M$ calling the resulting strand $4^\prime$.
The embedding of $M$ extends (by attaching disks to the $m_i$) to
an embedding of $B$ into the 3--ball containing the 4-component
string link. Thus, $4$ and $4^\prime$ cobound the singular capped
grope $\bar{g}^c$ from Construction~\ref{construct:sing-grope}
which sits inside $B$, where, by abuse of terminology, we let
$\bar{g}^c$ also denote the grope that has $4$ and $4^\prime$ as
its boundary.

Pick arcs $\alpha$ and $\beta$ contained in the bottom stage of
$\bar{g}^c$ and sharing endpoints with $4$ and $4^\prime$ such
that $\alpha\cup\beta$ splits $\bar{g}^c$ into three capped grope
cobordisms $g^c_1$, $g^c_2$ and $g^c_3$. If we number them
appropriately, $g^c_1$ modifies strand $4$ to the strand $\alpha$,
$g^c_2$ modifies $\alpha$ to $\beta$ and $g^c_3$ modifies $\beta$
to $4^\prime$. Note that each of these three capped grope
cobordisms is nonsingular.

Examining the way in which the caps hit the strands, we see that
$\amalg_j{\widetilde\tau}^c_3(G^c_j)=(+t_I)\amalg(-t_H)\amalg(+t_X)$,
where each $G^c_j$ is just $g^c_j$ together with trivial
cobordisms on the first three strands.

In order to get the desired $G^c$, we wish to glue these
cobordisms $G_i^c$ back together so that the resulting grope is
{\em embedded}. To do this, we use the transitivity argument from
\cite{CT1}, which is easily adapted to the current situation of
arcs rel boundary (as opposed to knots). In that argument the
individual gropes that are being glued together are homotoped
inside the ambient 3--manifold until they match up. However, the
homotopies are always isotopies when restricted to individual
gropes. (Except in the framing correction move where some twists
are introduced, which will not affect ${\widetilde\tau}^c(G^c)$.)
Thus ${\widetilde\tau}^c_3(G^c)=(+t_I)\amalg(-t_H)\amalg(+t_X)$ is
not changed during this procedure.

Now consider the case where $\ell=4$ but the univalent vertices of
the trees in the IHX relation are attached to strands $j_1$,
$j_2$, $j_3$ and $j_4$ which are not necessarily distinct. Then
the only modification needed in the above proof is to embed $M$ so
that the $m_i$ are meridians to the $j_i$th strand arranged in the
correct ordering ($i=1,2,3$), and make sure that the band from
$\partial\bar{g}$ attaches to the $j_4$th strand in the right
place.

Finally, if there are more than four strands, add the rest of the
strands to the picture away from the above construction.
\end{proof}

More generally, let us consider grope cobordisms of higher class.

\begin{thm}\label{ihxn}
Let $t_I$, $t_H$ and $t_X$ be three trees which differ by the
terms in an IHX relation in $\widehat{\cA}^t_n(\ell)$. Then there
is a class $n$ simple grope cobordism $G^c$, from the
$\ell$-component trivial string link to itself, such that
${\widetilde\tau}_n^c(G^c)=(+t_I)\amalg(-t_H)\amalg(+t_X)$.
\end{thm}
\begin{proof}
As in the proof of Theorem~\ref{thm:3dihx}, we will construct a
cobordism of one of the strands, extending the others by disks. As
argued at the end of the proof, it is sufficient to assume that no
two tips of any one tree are attached to the same component. Hence
we may assume that $\ell\geq n+1$. Further, as in that proof, we
may assume that $\ell=n+1$ on the nose.

Decompose $t_I$ into rooted trees $I,A,B,C,R$, where $I$
represents the ``I" in the IHX relation, a chosen root of $I$ is
connected to $R$, and the tips of $I$ connect to the roots of the
trees $A$, $B$ and $C$. Let the rooted tree given by $I$ union
$A$, $B$ and $C$ be called $t$ as illustrated in
Figure~\ref{treeI}.
\begin{figure}[ht!]
\centering
\includegraphics[width=2in]{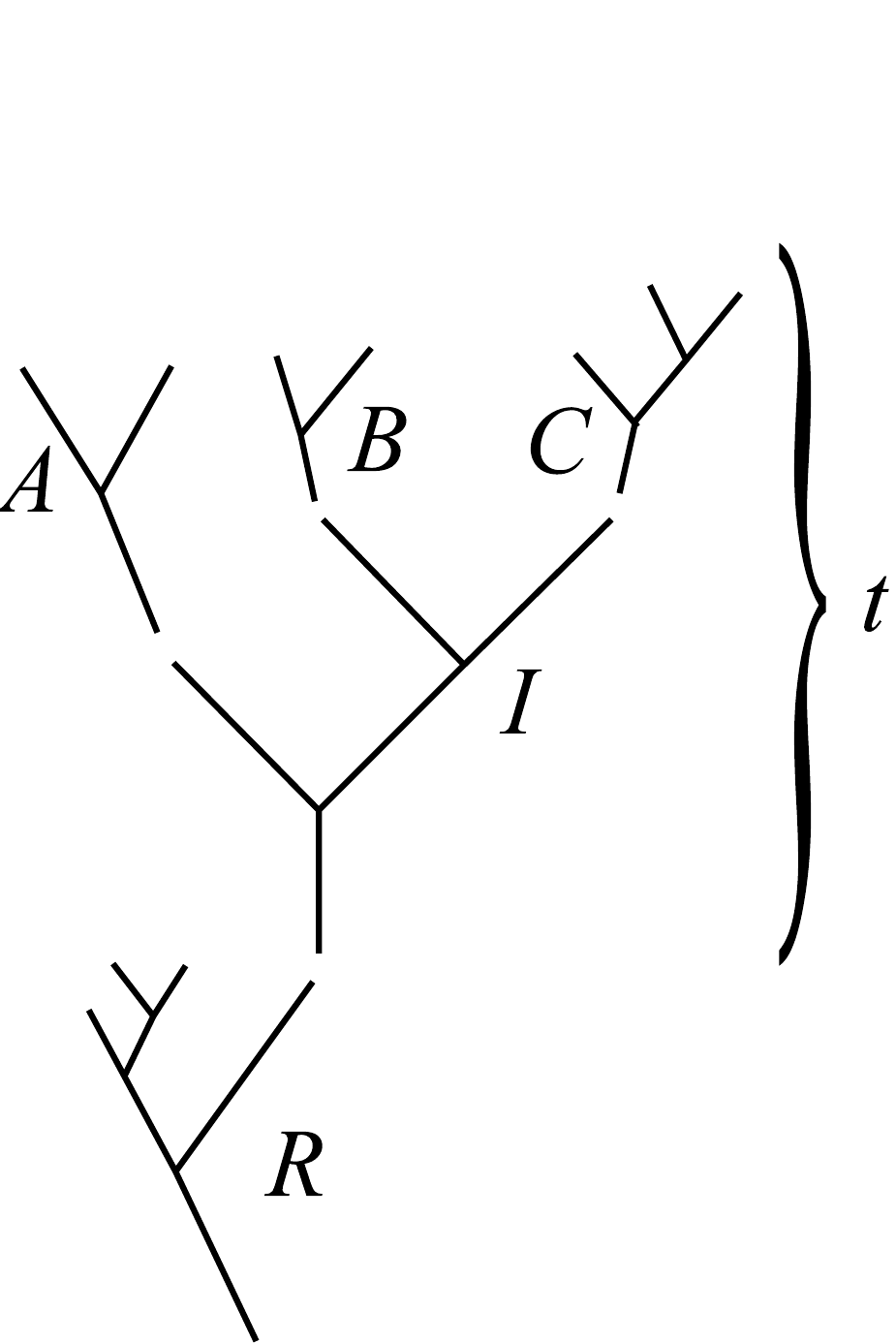}
\caption{}\label{treeI}
\end{figure}
Think of the ball containing $(n+1)$ strands as a
boundary-connected-sum $B_t\# B_R$, where $B_t$ is a ball with
strands which inherit the (distinct) labels of $t$ and $B_R$ is
ball with strands labeled  distinctly from the rest of
$\{1,2,\ldots, n+1\}$.

Consider a capped grope $g^c_t$ with one boundary component having
geometric intersection tree equal to the tree $t$ and contained in
$B_t$. (To see that such a grope exists note that a regular
neighborhood of a grope is a handlebody, which can be thought of
as a ball with a tubular neighborhood of some arcs removed. The
tips are part of a spine for the handlebody, so that there is a
bijection between tips and arcs, with each arc going through a
single tip once. Thus, the tips bound disks that are punctured by
distinct arcs. Now there is an embedding of this ball-with-arcs to
$B_t$ that takes the arcs to strands in $B_t$ according to any
bijection.)

Pruning the ``I" part $g^c_I$ of $g^c_t$, we get three capped
gropes realizing the trees $A$, $B$, $C$, denoted $g^c_A$,
$g^c_B$, $g^c_C$ respectively. As in Theorem~\ref{thm:3dihx},
consider the genus three handlebody $M$ which is the complement of
a trivial 3--strand string link with $m_i$ meridians to the
strands on $\partial M$. Taking $M$ to be a regular neighborhood
of $g_I^c$,
 there is an embedding of  $M$ into $B_t$ such
that the $m_i$ map to $\partial g^c_A$, $\partial g^c_B$,
$\partial g^c_C$. Now, by Construction~\ref{construct:sing-grope},
there is a singular grope $\bar{g}$ of class three inside $M$
which bounds an unknot on the boundary of $M$ such that the tips
of $\bar{g}$ bound parallel copies of $g^c_A$, $g^c_B$ and
$g^c_C$. (Note that these parallel copies intersect each other
because the second stages of $\bar{g}$ intersect each other, and
because parallel gropes in dimension three intersect.)

Let $g^c_R$ be a capped grope realizing the tree $R$ inside $B_R$,
such that the tip $T_0$ of $g^c_R$ corresponding to the tip of $t$
that connects to the root of $I$ bounds a cap that does not
intersect any strands. Note that $g^c_R$ can be surgered into a
disk, so that its boundary is unknotted.

Tube the cap on the (unknotted) tip $T_0$ of $g^c_R$ to the
(unknotted) $\partial \bar{g}$ on the boundary of $M$. Connect-sum
the (unknotted) $\partial g^c_R$ to (a push-off of) the strand in
$B_R$ corresponding to the root of $R$.

We get a singular capped grope cobordism $\bar{G}^c$ taking the
trivial $(n+1)$-component string link to itself. The connected
grope cobordism of the strand corresponding to the root of $R$ is
genus three at one stage and is embedded at that and all lower
stages (the ``$R$ part''). Higher stages (the ``$A$, $B$, and $C$
parts'') that lie above different genus one subsurfaces of the
genus three stage (in the ``$I$ part'') may intersect. Splitting
the grope via Proposition 16 of \cite{CT1}, we get three grope
cobordisms, each separately embedded, which can then be reglued by
transitivity, as in the proof of Theorem~\ref{thm:3dihx}, to get a
nonsingular grope cobordism, $G^c$, with
${\widetilde\tau}_n^c(G^c)=(+t_I)\amalg(-t_H)\amalg(+t_X)$.
\end{proof}

The previous theorem can be rephrased in the language of claspers
and implies Theorem~\ref{topihx} of the introduction.

A picture of three claspers of degree three as in
Theorem~\ref{topihx} is given in Figure~9 of \cite{CT2}. This was
derived from Theorem~\ref{thm:3dihx} using a mixture of claspers
and gropes in the following way. First, (using the notation in the
proof of Theorem~\ref{thm:3dihx}) the clasper representing $g^c_1$
was drawn. Next, we modified strand $1$ by $g^c_1$ to the new
position $\alpha$. We then drew in the clasper representing
$g^c_2$. This clasper intersects the grope $g^c_1$, but using the
usual pushing-down argument we pushed all the intersections down
to the bottom stage. We then pushed them off the strand $0$
boundary component of the grope, which is an isotopy in the
complement of $\alpha$. This gave rise to two disjoint claspers,
surgery on which moves strand $0$ to the arc $\beta$. The process
was repeated for the clasper representing $g^c_3$: it was pushed
out of the trace of the first two grope cobordisms/clasper
surgeries. We double-checked the result by performing surgery
along these three claspers and verified the result was isotopic to
the original trivial $4$-component string link.

\subsection{General IHX relations and the map $\widehat\tau^g_n$}
\label{sec:genihx} Next, we extend the realization of IHX
relations from trees to arbitrary diagrams. Extending the map
$\widehat\tau^c_n$ to {\em un}-capped grope cobordisms involves
some new wrinkles. First of all, in the absence of caps bounding
the grope tips, it will not be possible to attach the tips of the
grope-trees to $\ell$ strands with a meaningful ordering; however
tips will still be associated to components of the string link
according to the linking between the components and the
corresponding tips. Secondly, non-trivial linking between certain
tips will lead to the construction of graphs with non-zero Betti
number which result from gluing together the corresponding tips.

The reader may wonder why we do not introduce a map
$\widetilde{\tau}^g_n$ at the monoid level at this point. The
reason is that $\widetilde{\tau}^g_n$ is well-defined at the group
level, by Proposition~\ref{wdf:prop} below, but is not
well-defined at the monoid level, unless the choice of tips is
included as part of the grope data.

\begin{defi}\label{def:grope-degree-diagram-group}
Consider the abelian group generated by connected {\em diagrams}
(vertex-oriented unitrivalent graphs) whose univalent vertices are
labeled by the string link components $1,\cdots,\ell$ (possibly
with repeats), modulo the AS antisymmtery relations. Also divide
by the relation setting any diagram with a loop at a vertex to
zero.
 Let $\widehat
\cB^g_n(\ell)$ be the subgroup generated by such diagrams of grope
degree $n$. (Recall the grope degree is half the number of
vertices plus the first Betti number.)\hfill$\Box$
\end{defi}

\begin{rem} The fact that a loop at a vertex must be zero is a consequence
of IHX relations, provided that $n\geq 3$. In the case $n=2$, an
AS relation implies that such a diagram is $2$-torsion, and hence
is zero over any ring where $2$ is invertible.
\end{rem}

Now we define $\widehat\tau^g_n\colon\mathbb G_n(\ell)\ra
\widehat\cB^g_n(\ell)$. Let $G$ be a grope cobordism of class $n$.
First, choose a grope component $g\subset G$. As before, each
genus one branch of $g$ has an associated vertex-oriented
trivalent rooted tree $t$ whose tips $L_i$ correspond to tips
$T_i$ of $g$. For each such $T_i$, choose either a component $x_j$
of the string link, or another tip $T_j$ of $g$, and label the
corresponding tip $L_i$ of $t$ by $(L_i,x_j)$, or $(L_i,T_j)$
respectively. The root of $t$ is labeled  by the string link
component  that the boundary of $g$ meets. Now sum over all
choices to get a formal sum of labeled trees denoted $T(G)$.

Now we proceed to glue together some of the tips on each of these
labeled  trees, based on the geometric information of how the tips
link each other and the string link.
 Let $t$ be
a labeled  tree. It has tips $L_i$ labeled  $(L_i,T_j)$ or labeled
$(L_i,x_j)$, where each tip $L_k$ corresponds to the tip $T_k$. A
\emph{matching} of such a labeled  tree $t$ is a partition of the
set of all the tips of $t$ labeled by tips (and not string link
components) into pairs, such that the labels on each pair are of
the form $(L_i,T_j), (L_j,T_i)$. A matching determines a labeled
connected graph $\Gamma$, gotten by gluing together matched tips
of $t$, where each edge resulting from such a gluing assumes the
coefficient $\lk (T_i,T_j)=\lk (T_j,T_i)$. Each of the remaining
univalent vertices $L_i$ is labeled by some component $x_j$, and
assumes the coefficient $\lk(T_i,x_j)$. Each such $\Gamma$
determines a multiple of a generating diagram of
$\widehat\cB^g_n(\ell)$, where the coefficient of the diagram is
the product of the coefficients on the tips and edges of $\Gamma$.
Define $\langle t \rangle$ as the sum of these elements in
$\widehat\cB^g_n(\ell)$ over all matchings of $t$. If there are no
matchings, then $\langle t \rangle \,=0$ by definition. Extend
$\langle \cdot \rangle$ to linear combinations of trees linearly.
Now define $\widehat\tau^g_n(G)$ to be $\langle T(G)\rangle$.

\begin{rem}
\noindent\begin{enumerate}
\item If $G$ extends to a simple grope cobordism $G\subset
G^c$, then $\widehat\tau^g_n(G)$ is just the image of
$\widehat\tau^c_n(G^c)$ under the map $\poff\colon
\widehat\cA_n^t(\ell)\ra \widehat\cB_n^g(\ell)$ that pulls the
trees off the components of $\ell$ and labels their univalent
vertices accordingly.
\item If one translates a grope cobordism into a union of rooted
clasper surgeries, the map $\widehat\tau^g_n$ can be calculated as
follows. Instead of $T(G)$, consider the associated tree of each
clasper with root labeled by the strand linked by the clasper's
root, and then apply $\langle\cdot\rangle$ as before. Then sum
over all of the claspers. If the rooted clasper, $C$, can be
turned into a simple clasper, $C'$, by turning Hopf pairs of tips
into edges, then $\widehat\tau^g_n(C)$ is the diagram which is the
associated graph of  $C'$, with univalent vertices labeled
according to where the capped tips of $C'$ meet the string link.
\hfill$\Box$
\end{enumerate}
\end{rem}

\begin{prop}\label{wdf:prop}
The map $\widehat\tau^g_n$ is well-defined.
\end{prop}

We prove this in the appendix, where we consider the more general
situation of gropes which may not be of genus one.

Finally, we show that the IHX relation can be realized in the
world of graphs by uncapped gropes.

\begin{thm}\label{genihx}
Let $D_I,D_H,D_X\in \widehat\cB _n^g(\ell)$ be diagrams differing
by the terms in an IHX relation. Then there is a grope cobordism
$G$, from the trivial $\ell$-string link to itself, such that
$\widehat\tau^g_n(G)=D_I-D_H+D_X$.
 \end{thm}
\begin{proof}  First, cut
some edges (not contained in the ``I'' part) of $D_I$ to make a
tree $D_I^t$. Pick a univalent vertex that did not come from a cut
as the root. Let $\ell$ be the number of tips. As before, think of
the complement of a trivial $\ell$ string link as a handlebody,
$M$, with special curves $\{m_i\}_{i=1}^{\ell}$ on its boundary.
Let the tips of $D_I^t$ be placed in correspondence with the
curves $m_i$. Embed $M$ in the complement of a trivial string
link, such that if a tip $L_i$ of $D_I^t$ is labeled  by a
component $x$ of the string link, then the corresponding $m_i$
links $x$ exactly once. Also, tips resulting from cuts of $D_I$
should have the corresponding $m_i$  linking exactly once. Take a
trivial subarc of the component of the string link corresponding
to the root of $D_I^t$ and perform a finger move so that it goes
through $M$ as a trivial subarc $\eta$. Now the proof of
Theorem~\ref{ihxn} yields a ``weak'' capped grope cobordism
$g^{\bar{c}}$ (with $g\subset M$) which modifies $\eta$, where the
weakness comes from the fact that here the linking pairs of tips
have intersecting caps. Ignoring this defect, $g^{\bar{c}}$
extends (as in the proof of Theorem~\ref{ihxn}) to a (weak) capped
grope cobordism $G^{\bar{c}}$ of the trivial string link such that
$\widehat\tau^g_n(G^{\bar{c}})= D_I-D_H+D_X$. This can be seen as
follows. Note that in this case $\widehat\tau^g$ behaves just like
$\widehat\tau^c$, except that it identifies tips corresponding to
Hopf-linked tips (where the caps intersect), and hence glues the
cut edges back together. The three different genus one pieces of
$G^{\bar{c}}$ link with each other in rather a complicated way but
this is not seen by the map $\widehat\tau^g_n$. Also note that the
tips of $G$ are parallel to the curves $m_i$, so that the map
$\widehat{\tau}^g_n$  labels the univalent vertices appropriately.
\end{proof}

\section{Grope cobordism of string links}\label{sec:string}

Let  $\cL(\ell)$ be the set of isotopy classes of string links in
$D^3$ with $\ell$ components (which is a monoid with respect to
the usual ``stacking'' operation). Its quotient by the relation of
grope cobordism (respectively capped grope cobordism) of class $n$
is denoted $\cL(\ell)/G_n$ (respectively $\cL(\ell)/G^{c}_n$),
compare Definition~\ref{grope cobordism}. The submonoid of
$\cL(\ell)$, consisting of those string links which cobound a
class $n$  grope (respectively capped grope) with the trivial
string link, is denoted by $G_n(\ell)$ (respectively
$G^{c}_n(\ell)$).

\begin{proof}[Proof of Theorem~\ref{thm:nil}]
Let us begin with the statements for the capped case. Then
$\cL(\ell)/G^{c}_{n}$ can be identified with the quotient of
$\cL(\ell)$ modulo the relation of simple clasper surgery of
class~$n$. This translation works just like for knots where it was
explained in \cite{CT1}. All the results then follow from
\cite[Thm.5.4]{H}. For example, the fact that the iterated
quotients are central is proven by showing that $ab=ba$, modulo
simple clasper surgery of class~$(n+1)$, if $a$ is a string link
that is simple clasper n-equivalent to the trivial string link.
This follows by sliding the claspers (that turn the trivial string
link into $a$) past another string link $b$.

In the absence of caps one has to translate into rooted clasper
surgery of grope degree~$n$ instead, as explained in \cite{CT1}.
Just as above, all results follow from the techniques of Habiro
\cite{H}.
\end{proof}

This result makes it possible to try to compute the {\em abelian}
iterated quotients in terms of diagrams, which we proceed to do.
We shall first define the map from diagrams to string links modulo
grope cobordism:
$$
\Phi_n\colon\cB^g_n(\ell)\to G_n(\ell)/G_{n+1}.
$$

Indeed, we defined this for $\ell = 1$ in \cite{CT2} in the
following way. Given a diagram $D\in \widehat\cB^g_n(\ell)$, find
a grope cobordism $G$ of class~$n$, \emph{corresponding to a
simple clasper}, such that $\widehat\tau^g_n(G) = D$. Then define
\[
\widehat\Phi_n(D)=\partial_1G(\partial_0G)^{-1},
\]
 where $\partial G =\partial_0G\cup \partial_1G$. One must show that the map is
well-defined, i.e. that the choice of embedding of the simple
clasper does not matter. The argument given in \cite{CT2} works
 with little modification for all $\ell \geq 1$.

The next proposition implies that we can take any grope cobordism
$G$ satisfying $\widehat\tau^g_n(G)=D$ in the above definition,
not having to restrict to those corresponding to simple claspers.

\begin{prop}\label{wdf}
 Given any grope cobordism $G$ of class~$n$,
 $\partial_1G(\partial_0G)^{-1}=\widehat\Phi_n\circ \widehat\tau^g_n(G)\in G_n(\ell)/G_{n+1}$
\end{prop}
\begin{proof}

 Any grope cobordism can be refined to a sequence of genus one grope cobordisms by
Proposition~16 of \cite{CT1} and this refinement evidently
commutes with $\widehat\tau^g_n$. Then, using Theorem~35 of
\cite{CT1}, each of these cobordisms can be refined into a
sequence of simple clasper surgeries and clasper surgeries of
higher degree, and this refinement commutes with
$\widehat{\tau}^g_n$. (To see this it suffices to notice that the
``zip construction" commutes with $\widehat{\tau}^g_n$.) Thus
$$
\partial_1G(\partial_0G)^{-1} = (\partial_1G)(L_k)^{-1}(L_k)(L_{k-1})^{-1}\cdots...
(L_1)(\partial_0G)^{-1},
$$
where the $L_i$ are string links modified by successive simple
clasper surgeries. Note that we can omit any pairs
$(L_i)(L_{i-1})^{-1}$ corresponding to clasper surgeries of higher
degree, since this product is trivial in $\cL(\ell)/G_{n+1}$. On
the other hand, we know that for pairs  differing by simple
claspers $C_i$ of degree $n$, $(L_i)(L_{i-1})^{-1}=
\widehat\Phi_n(\widehat\tau^g_n(C_i))$, by definition of
$\widehat\Phi_n$. Thus
\begin{align*}
\partial_1G(\partial_0G)^{-1}  &=  \#_i \widehat\Phi_n(\widehat\tau^g_n(C_i))\\ &=\widehat\Phi_n(\widehat\tau^g_n(\sum
C_i))\\ &=\widehat\Phi_n(\widehat\tau^g_n(G))
\end{align*}
which completes the proof.
\end{proof}

We next show that $\widehat\Phi_n$ vanishes on all IHX relations
and hence descends to a well-defined map $\Phi_n$.

\begin{thm}\label{homom}
$\Phi_n\colon \cB^g_n(\ell)\to G_n(\ell)/G_{n+1}$ is a
well-defined surjective  homomorphism.
\end{thm}

\begin{proof}[Proof of Theorem \ref{homom}] By Theorem~\ref{genihx}, any IHX relation, $R_{IHX}$,
is the image under $\widehat\tau^g_n$ of a grope cobordism, $G$,
from a trivial string link to another trivial string link, denoted
$1_\ell$. So by Proposition~\ref{wdf},
\begin{align*}
\widehat\Phi_n(R_{IHX})&=\widehat \Phi_n(\widehat\tau^g_n(G))\\
&=(\partial_1G)(\partial_0G)^{-1}\\ &=1_\ell\# 1_\ell^{-1}\\
&=1_\ell
\end{align*}
Next we consider surjectivity of $\Phi_n$. The elements of
$G_n(\ell)$ are by definition of the form $\partial_1 G$ where $G$
is a class $n$ grope cobordism with $\partial_0G=1_\ell$. By
Proposition~\ref{wdf}, $\partial_1G
=\Phi_n\circ\widehat\tau^g_n(G)$.

\end{proof}

Using the Kontsevich integral as a rational inverse, we are now
able to prove Theorem~\ref{kontsevich-grope} which says that
$\Phi_n$ turns into an isomorphism after tensoring with $\Q$.

\begin{proof}[Sketch of proof of Theorem~\ref{kontsevich-grope}]
This was proven in full detail in \cite{CT2} for the case when
$\ell=1$. One sets up the (logarithm of the) Kontsevich integral
as an inverse. Using the Aarhus integral \cite{A}, it is easy to
show that the bottom degree term of the Kontsevich integral
coincides with our map $\widehat\tau^g_n$. More precisely, if $G$
is a grope cobordism, then Aarhus surgery formulae show that
$$(\log Z)_n(\partial_1G(\partial_0G)^{-1})=\widehat\tau^g_n(G),$$
where $(\log Z)_n$ is of grope degree $n$. Thus $\Phi_n((\log
Z)_n(\partial_1G(\partial_0G)^{-1}))=\partial_1G(\partial_0G)^{-1}$,
or $\Phi_n\circ (\log Z)_n = id$. On the other hand $(\log
Z)_n(\Phi_n(D))=(\log Z)_n(\partial_1G(\partial_0G)^{-1})$ for a
grope $G$ satisfying $\widehat\tau_n(G)=D$. But then, by the above
highlighted formula we can conclude that $(\log
Z)_n\circ\Phi_n=id$.

Also, the Kontsevich integral of grope cobordisms of class $(n+1)$
will lie in degree $(n+1)$, so that the Kontsevich integral indeed
factors through $G_n(\ell )/G_{n+1}\otimes \Q$. (Here we use the
fact that the Kontsevich integral of string links preserves the
loop (and hence grope) degree.) The fact that $\log Z_n$ is a
homomorphism is straightforward using the Aarhus formula. (In
\cite{CT2} we used the Wheeling isomorphism to show this for
knots, but that was unnecessary. The lowest degree part of the
Wheeling isomorphism is just the identity.)
\end{proof}

It is unknown whether the analogous statements for the relation of
{\em capped} grope cobordism of string links are true. There are
two difficulties, one is the question of whether one can realize
the STU-relations in $\cA_n(\ell)$ by capped grope cobordisms. The
other is the question of whether Habiro's main theorem \cite{H}
generalizes from knots to string links: Does the Vassiliev
filtration of string links agree with the relation generated by
simple clasper surgery? It follows from the techniques of
\cite{CT1} that the latter agrees with capped grope cobordism.

We conclude this section by carefully proving Lemma 3.11 (c) from
\cite{CT2}, which we restate here for convenience.
\begin{lem}\label{ct2lem}
Let $U$ be the unknot. Suppose three claspers $C_i$ of grope
degree $n$ on $U$ differ according to the IHX relation. Then
$U_{C_1}\# U_{C_2}\# U_{C_3}\in G_{n+1}(1)$.
\end{lem}
\begin{proof}
Let $K=U_{C_1}\# U_{C_2}\# U_{C_3}$. The union of the three
claspers corresponds to a grope cobordism, $g$, of class $n$
between the unknot and $K$,  where the bottom stage is of genus
three. By Proposition
 \ref{wdf}, we have that $K={\Phi}_n\circ \widehat{\tau}^g_n(g)$. However
 $\widehat{\tau}^g_n(g)$ is an IHX relator, and so by Theorem \ref{homom},
 ${\Phi}_n$ vanishes on it. Thus $K$ is trivial in $G_n(1)/G_{n+1}$, implying that
 $K\in G_{n+1}(1).$
\end{proof}

\section{Appendix: Associating a linear combination of graphs to an arbitrary grope}
In this appendix we consider the set of class $n$ grope cobordisms
of $\ell$-string links, which may not be of genus one. Let this
set be denoted $\widehat{\mathbb G}_n(\ell)$.

Now we define $\widehat\tau^g_n\colon \widehat {\mathbb
G}_n(\ell)\ra \widehat\cB^g_n(\ell)$. Let $G$ be a grope cobordism
of class $n$. First, choose a grope component $g\subset G$. Choose
tips for the grope component. Associate a linear combination of
trees to $g$ as follows. Each stage of $g$ has a set of hyperbolic
pairs of basis elements which bound further stages of the grope,
or are tips. A branch of the grope is defined to be a choice of
such a pair at the bottom stage, followed by a choice of
hyperbolic pair at each stage which is bounded by the first pair,
and so on. Each branch of the grope has an evident tree associated
with it, whose tips $L_i$ correspond to the tips $T_i$ of the
branch of the grope. As in the construction after
definition~\ref{def:grope-degree-diagram-group} in
Section~\ref{sec:genihx}, for each such $T_i$, choose either a
component $x_j$ of the string link, or another tip $T_j$ of $g$,
and label the corresponding tip $L_i$ of $t$ by $(L_i,x_j)$, or
$(L_i,T_j)$ respectively. The root of $t$ is labeled  by the
component of the string link that the boundary of $g$ meets. Now
sum over all choices, including all choices of branch of $g$, to
get a formal sum of labeled trees denoted $T(G)$.

Now define $\widehat\tau^g_n(G)$ to be $\langle T(G)\rangle$, as
before.

\begin{prop}\label{wdf:prop}
$\widehat\tau^g_n$ is well-defined.
\end{prop}
\begin{proof}
The first ambiguity is the orientation. As in Lemma~\ref{lem:AS
for gropes}, changing a positive quadrant results in a change of
orientation all of the higher stages, including pushing annuli.
Changing the orientation of a pushing annulus changes the sign of
every term in $\widehat\tau^g_n$, either by reversing the sign of
the linking number with another tip, or by changing the sign of
the linking number with a string link component. Thus, as in the
proof of Lemma~\ref{lem:AS for gropes}, there are an even number
of sign changes.

The second ambiguity is the choice of pushing annuli. Every tamely
embedded grope can be extended to include pushing annuli, but this
extension may not be unique. At a top stage of the grope, there
will be two choices for every hyperbolic pair of tips,
 according to whether a given annulus extends ``up" or ``down" from the surface stage. Changing the choice of pushing annuli at a hyperbolic pair of tips has the effect of switching which quadrants are positive. However, the cyclic order of the vertex does not change. The induced orientations of the pushing annuli are either the same, or they are both reversed, resulting in no net change in sign.

The third ambiguity arises from choosing different tips for a
grope component $g\subset G$. Notice that $\widehat\tau^g_n$ never
sees the linking of tips on the same stage of $g$: Either they
belong to different branches and hence will be part of different
tree summands, or they are dual to each other in which case a
graph with a loop at a vertex would result. Thus on a single
surface stage, the linking number with objects $c_i$ is all that
matters, where $c_i$ is either a component of the string link or
another tip of $g$ on a different stage.

Suppose we are not at a top stage. Then at least one curve in
every hyperbolic pair bounds a higher surface stage. Removing a
regular neighborhood of the higher surface stages, we get a planar
surface. The tips become arcs joining some pairs of boundary
components. Different choices of tips are related by Dehn twists
on curves in the planar surface. Note that the boundary components
of the planar surface are all null-homologous in the complement of
$\cup c_i$. (They bound surfaces, and if the surfaces are slightly
perturbed, they avoid $c_i$.) Hence choices of tips differ by
multiples of curves which link the $c_i$ trivially and hence do
not change the contribution of $g$ to $\widehat\tau^g_n(G)$.

Now suppose we are at a top stage of genus $m$. Any two choices of
tips = symplectic bases $(\alpha_1,
\beta_1,\ldots,\alpha_m,\beta_m)$ are related by an element of
$Sp(2m, \mathbb Z)$, which is generated by the following
automorphisms:
\begin{itemize}
\item
for some $i$, $\alpha_i\mapsto \alpha_i+\beta_i $ and everything
else is fixed
\item for some $i$, $\beta_i\mapsto \alpha_i+\beta_i $ and everything else
is fixed
\item for some $i\neq j$, $\begin{cases} \alpha_i\mapsto \alpha_i+\alpha_j \\
\beta_j\mapsto -\beta_i+\beta_j\end{cases} $ and everything else
is fixed
\item  for some $i\neq j$ $\begin{cases}\alpha_i\mapsto \alpha_i+\beta_j \\
\alpha_j\mapsto \beta_i+\alpha_j \end{cases}$ and everything else
is fixed
\item  for some $i\neq j$$\begin{cases}\beta_i\mapsto \beta_i+\alpha_j \\
\beta_j\mapsto \alpha_i+\beta_j \end{cases}$ and everything else
is fixed
\item  for some $i\neq j$ $\begin{cases}\beta_i\mapsto \beta_i+\beta_j \\
\beta_j\mapsto -\alpha_i+\alpha_j \end{cases}$ and everything else
is fixed
\end{itemize}

Let us adopt the following notation for expressing the
contribution $T(g)$ of $g$ to $T(G)$. Compute the disjoint union
of trees where the tips correspond to the tips of $g$, and label
each tip $L_i$ by a linear combination $\sum_r n_r \mathbf{c_r}$
where the labels $\mathbf{c_r}$ correspond to components of the
string link and tips $T_j$ of $g$ with $j\neq i$ (and $n_r$ is the
corresponding linking number with $T_i$). This represents $T(g)$
by expanding the trees linearly in the labels. Note that if any
labeled  trees in $T(g)$ represent zero modulo AS or IHX
relations, then these relations will still be present upon gluing,
so that the corresponding contribution to
$\widehat\tau^g_n(G)=\langle T(G)\rangle$ will also be zero.

The trees in $T(g)$ before and after applying the first
automorphism above only differ in a subtree isomorphic to a ``Y",
which we can represent by a bracket $[\, , \, ]$. The difference
is then represented by
$$
\left[\sum_r \lk(\alpha_i,c_r) \mathbf{c_r},\sum_r
\lk(\beta_i,c_r)\mathbf{c_r}\right] - \left[\sum_r
\lk(\alpha_i+\beta_i,c_r)\mathbf{c_r},\sum_r \lk(\beta_i,
c_r)\mathbf{c_r}\right].
$$
Breaking the second summand into two terms, we see that that
$$
\left[\sum_r \lk(\beta_i,c_r)\mathbf{c_r},\sum_r
\lk(\beta_i,c_r)\mathbf{c_r}\right]=0,
$$
 is sufficient to show that $T(g)$, and hence
$\widehat\tau^g_n(G)$, remains unchanged. The fact that $[x,x]=0$
corresponds to the statement that a loop at a vertex is zero. The
case of the second automorphism is handled in the same way.

Let's consider the third automorphism. Abbreviate the notations
$\sum_r \lk(\alpha, c_r)\mathbf{c_r}$ by $\lk(\alpha,c)$. Then
notice that the difference in $T(g)$ only occurs in the $i$ and
$j$ trees, and this difference is
\begin{align*} \left[\lk(\alpha_i,c),\lk(\beta_i,c)\right]+\left[\lk(\alpha_j,c),\lk(\beta_j,c)\right] -
\left[\lk(\alpha_i+\alpha_j,c),\lk(\beta_i,c)\right]- \phantom{xxxxxxxxxxxxxx}\\
\left[\lk(\alpha_j,c),\lk(-\beta_i+\beta_j,c)\right], \end{align*}
which is easily seen to be zero. The cases of the last three
automorphisms are handled identically.
\end{proof}

We remark that Proposition~\ref{wdf} is still true for this
extended definition of $\widehat\tau^g_n$.

\end{document}